\documentclass[12pt,a4paper]{article}

\usepackage[utf8]{inputenc}
\usepackage{amsmath}
\usepackage{tikz-cd}
\usepackage[all]{xy}
\usepackage{amssymb}
\usepackage{amsthm}
\usepackage{float}
\usepackage{multicol}
\usepackage{graphicx}
\usepackage{color}

\usepackage{hyperref}

\newtheorem{theorem}{Theorem}[section]
\newtheorem*{theorem*}{Theorem}

\newtheorem{lemma}[theorem]{Lemma}
\newtheorem{proposition}[theorem]{Proposition}
\newtheorem{definition}[theorem]{Definition}
\theoremstyle{definition}
\newtheorem{remark}[theorem]{Remark}
\newtheorem{notation}[theorem]{Notation}
\newtheorem{example}[theorem]{Example}
\newtheorem{conjecture}[theorem]{Conjecture}

\title{From Differential Values to Roots of the Bernstein-Sato Polynomial}
\author{David Senovilla-Sanz\footnote{The author is supported by the Spanish research project PID2022-139631NB-I00 funded by the Agencia Estatal de Investigación - Ministerio de Ciencia e Innovación. Moreover, he is also supported by a predoctoral contract “Concepción Arenal” of the  Universidad de Cantabria}}

\newcommand{\Addresses}{{
  \bigskip
  \footnotesize

  David Senovilla-Sanz. \textsc{Departamento de Matemáticas, Estadística y Computación. Universidad de Cantabria. Avda. de los Castros s/n, 39005 -- Santander, SPAIN}\par\nopagebreak
  \textit{E-mail address}, D. Senovilla-Sanz: \texttt{david.senovilla@unican.es}

}}

\begin{document}
\maketitle
\begin{abstract}

Let $C$ be a cusp in $(\mathbb C^2,\mathbf 0)$ with Puiseux pair $(n,m)$. This paper is devoted to show how the semimodule of differential values of $C$ determines a subset of the roots of the Bernstein-Sato polynomial of $C$. We add more precise results when
the multiplicity of the cusp is $n\leq 4$.
\end{abstract}

\tableofcontents

\section{Introduction}

In this paper we study the relationship between two analytic invariants associated with a plane branch $C$ in $(\mathbb C^2,\mathbf 0)$: the semimodule of differential values of $C$ and the Bernstein-Sato polynomial of $C$. Specifically, we show how the semimodule of differential values of $C$ determines a subset of the roots of the Bernstein-Sato polynomial of $C$, when $C$ is a cusp. We studied this relationship motivated by following two facts.

First, the semimodule of differential values of a branch appears as a crucial element in the analytic classification of plane branches, as described in \cite{hefez2}. More precisely, determining whether two branches are analytically equivalent requires, as the first step, to compute their respective semimodules of differential values. 

Second, the roots of the Bernstein-Sato polynomial are directly related with the eigenvalues of the monodromy map, see \cite{veys2024}. This shows its relevance when studying the monodromy conjecture. Moreover, \cite{blanco} shows that, in the case of curves, the Bernstein-Sato polynomial determines the poles of the complex zeta function of $C$.

\bigskip

To clarify these results, we begin by recalling the definition of the semimodule of differential values.

Let $C$ be a germ of a holomorphic branch in $(\mathbb C^2,\mathbf 0)$. Consider a Puiseux parametrization $\phi(t)=(x(t),y(t))$ of $C$, where $x(t),y(t)\in\mathbb C\{t\}$ are one variable convergent complex power series, with the extra condition that the multiplicity of $C$ at the origin $\mathbf 0$ is  
$$
\nu_\mathbf 0(C)= \min\{\operatorname{ord}_t(x(t)), \operatorname{ord}_t(y(t))\}.
$$ 
Here $\operatorname{ord}_t$ denotes the order of a series in the variable $t$.

Given a function $h\in\mathbb C\{x,y\}$, we denote by $\nu_C(h):=\operatorname{ord}_t(h\circ \phi)$. The semigroup of $C$ is then defined as
$$
\Gamma_C:=\{\nu_C(h): h\in\mathbb C\{x,y\}\}.
$$
Similarly, given a holomorphic 1-form $\omega\in\Omega_{\mathbb C^2,\mathbf 0}^1$, we write the pull-back of $\omega$ by $\phi$ as $\phi^*\omega=\alpha(t)dt$. Again, we denote by $\nu_C(\omega):=\operatorname{ord}_t(\alpha)+1$ the differential value of $\omega$. The semimodule of differential values of $C$ is then
$$
\Lambda_C:=\{\nu_C(\omega):\omega\in\Omega_{\mathbb C^2,\mathbf 0}^1\}.
$$
It is always true that the set $\Lambda_C\setminus \Gamma_C$ is finite.
\bigskip 

The Bernstein-Sato polynomial of the branch $C$ is defined as follows. Let $f\in \mathbb C\{x,y\}$ be an implicit equation of $C$, and denote by $\mathcal D$ the non commutative ring of differential operators in variables $x,y$. Introducing a new variable $\rho$, the action of $\mathcal D$ over $\mathbb C\{x,y\}$ extends in a natural way to elements of the of form $h^\rho$, where $h\in\mathbb C\{x,y\}$. The set of polynomials $B\in \mathbb C[\rho]$ for which there exists an operator $P\in\mathcal D[\rho]$ satisfying 
$$
P(\rho)f^{\rho+1}=B(\rho)f^\rho,
$$
is a non zero principal ideal $I\subset \mathbb C[\rho]$, see \cite{bjork}. The monic generator of this ideal is called the Bernstein-Sato polynomial of $C$.

Assume that $C$ is a cusp, in other words, the semigroup of $C$ is given by $\Gamma_C=\langle n,m\rangle$, with $gcd(n,m)=1$ and $2\leq n<m$. Denote by $\lambda_1=\min(\Lambda_C\setminus \Gamma_C)$, in case it exists. This number $\lambda_1$ is equivalent to the Zariski's invariant of $C$ from \cite{zariski}. We show the following.
\begin{theorem}\label{thm:bern:l1}
Let $C$ be a cusp with semigroup $\Gamma_C=\langle n,m\rangle$ and semimodule of differential values $\Lambda_C$. Assume that $\lambda_1=\min(\Lambda_C\setminus\Gamma_C)$ exists. Then for any element $\lambda\in (\lambda_1+\Gamma_C)\setminus \Gamma_C\subset \Lambda_C$, the rational number $-\lambda/nm$ is a root of the Bernstein-Sato polynomial of $C$.
\end{theorem}
\begin{theorem}\label{thm:bern:4}
Let $C$ be a cusp with semigroup $\Gamma_C=\langle n,m\rangle$ and semimodule of differential values $\Lambda_C$. Assume that $n\leq 4$. Then for any element $\lambda\in \Lambda_C\setminus \Gamma_C$, the rational number $-\lambda/nm$ is a root of the Bernstein-Sato polynomial of $C$.
\end{theorem}

Based on these results, we conjecture that Theorem \ref{thm:bern:4} holds for any $n$. Specifically:

\begin{conjecture}
Let $C$ be a cusp with semigroup $\Gamma_C=\langle n,m\rangle$ and semimodule of differential values $\Lambda_C$. Then for any element $\lambda\in \Lambda_C\setminus \Gamma_C$, the rational number $-\lambda/nm$ is a root of the Bernstein-Sato polynomial of $C$.
\end{conjecture}

The proofs of both theorems rely on the algebraic conditions given in \cite{cassou} to obtain roots of the Bernstein-Sato polynomial of $C$. These algebraic conditions are compared with those derived when requiring a given natural number to belong to $\Lambda_C$. These last conditions are determined using the classical Büchberger's algorithm for computing minimal standard bases of ideals in $\mathbb C\{x,y\}$. 

Now we proceed to recall the definitions of the last concepts.

\section{Preliminaries}
\subsection{Standard basis for an ideal}\label{sec:ideal}
In this first subsection we recall the main concepts used about standard basis. We follow \cite{hefez1} along all the subsection. Let $p$  be a positive integer, we consider the semigroup $(\mathbb Z_{\geq 0})^p$, where we denote by $\mathbf 0=(0,\ldots,0)$ the neutral element in $(\mathbb Z_{\geq 0})^p$. Let $\preceq$ be a \emph{monomial order} of $(\mathbb Z_{\geq 0})^p$, that is, $\preceq$ is a total order in $(\mathbb Z_{\geq 0})^p$ satisfying the following conditions:
\begin{itemize}
\item $\mathbf 0\preceq s$ for all $s\in(\mathbb Z_{\geq 0})^p$.
\item given $s_1\preceq s_2$, then $s_1+s\preceq s_2+s$ for all $s,s_1,s_2\in (\mathbb Z_{\geq 0})^p$.
\end{itemize}
We also ask $\preceq$ to satisfy the for all $s_0\in(\mathbb Z_{\geq 0})^p$, the set $\{s\preceq s_0:s\in(\mathbb Z_{\geq 0})^p\}$ is finite. In other words, we demand $\preceq$ to have the \emph{finiteness property}.

\begin{example}\label{ex:pr:orders}
The following monomial order is the main one that we are going to use along this paper. 
Consider a pair of non negative integers $(n,m)$, with $2\leq n<m$ and without common factor. Then given $(a,b),(c,d)\in(\mathbb Z_{\geq 0})^2$, we say that $(a,b)\prec (c,d)$ if and only if either $na+mb<nc+md$ or $na+mb=nc+md$ with $a<c$. We call this order the \emph{weighted order} with respect the weights $(n,m)$.
\end{example}

Given $s,t\in (\mathbb Z_{\geq 0})^p$, we say that  $s$ \emph{divides} $t$ or that $t$ is divisible by $s$, if $t-s\in(\mathbb Z_{\geq 0})^p$.

We denote by $\mathbb C\{x_1,\ldots,x_p\}$ the ring of convergent power series with $p$ variables. When $p=2$ we write $\mathbb C\{x,y\}$ instead of $\mathbb C\{x_1,\ldots,x_p\}$. For an element $g\in \mathbb C\{x_1,\ldots,x_p\}$, we write 
$$
g=\sum_{\alpha\in I}a_\alpha x^\alpha, \quad I\subset (\mathbb Z_{\geq 0})^p, \quad x^\alpha=x_1^{\alpha_1}x_2^{\alpha_2}\ldots x_n^{\alpha_n};
$$ 
we denote the \emph{Newton cloud} of g by
$$
\mathcal N(g):=\{\alpha\in(\mathbb Z_{\geq 0})^p: a_\alpha\neq 0\}.
$$
If $g\neq 0$, then we define the \emph{leading power} of $g$ with respect to $\preceq$ as $lp(g):=\min(\mathcal N(g))$. In the case where $g=0$, we set $lp(0)=(\infty,\infty,\ldots,\infty)\succ \alpha$ for any $\alpha\in(\mathbb Z_{\geq 0})^p$. Now, the \emph{leading term} of $g$ is $lt(g):=a_\beta x^\beta$ with $\beta=lp(g)$ if $g\neq 0$ and $lt(0)=0$.
\begin{definition}
Let $I\subset \mathbb C\{x_1,\ldots,x_p\}$ be an ideal and $B\subset I$ a subset. We say that $B$ is a \emph{standard basis} of $I$ if:
\begin{itemize}
\item $B$ generates $I$ as an ideal.
\item For any $h\in I$, there exists $b\in B$, such that $lp(b)\mid lp(h)$.
\end{itemize}
We say that $B$ is a \emph{minimal standard basis}, if for all $b\in B$, then $B\setminus\{b\}$ is not a standard basis. 
\end{definition}
 
We now recall some terminology necessary to understand the Büchberger's algorithm for computing minimal standard bases. Let $B\subset \mathbb C\{x_1,\ldots,x_p\}$ be subset, and consider two elements $g,r\in \mathbb C\{x_1,\ldots,x_p\}$. We say that $r$ is a \emph{reduction} modulo $B$ of $g$, if there exist $a\in\mathbb C$, $\alpha\in(\mathbb Z_{\geq 0})^p$ and $b\in B$, such that
$$
r=g-ax^\alpha b,
$$
where either $r=0$ or $lp(g)\succ lp(r)$. If such an $r$ exists, we say that $g$ is reducible modulo $B$.

We denote by $r_\infty$ a \emph{final reduction} of $g$ as the Krull limit of a sequence of reductions starting of $g$ until one of the following situations occurs: either the zero element is reached, or we obtain an element no longer non reducible by $B$. 

An element $r'\in \mathbb C\{x_1,\ldots,x_p\}$ is called a \emph{partial reduction} of $g$ modulo $B$, if there exists a finite sequence of reductions starting at $g$ that gives $r'$. 

In both the final and partial reduction definitions, we include the case where $g$ is non reducible modulo $B$; in this situation, $g$ itself is considered its own final (resp. partial) reduction.

As an abuse of notation, we also say that $r$ is a final reduction of $g$ modulo $(B)$, if $r=g-h$ with $h$ belongs to the ideal $(B)$ and $r$ is non reducible modulo $B$.

\begin{remark}\label{rem:pr:pertenencia ideal}
Assume that $B$ is a finite set. By Artin's Approximation Theorem, see \cite{artin}, saying that a final reduction of $g$ modulo $B$ is zero is equivalent to saying that $g$ belongs to the ideal generated by $B$.
\end{remark}
Take non zero functions $g_1,g_2\in \mathbb C\{x_1,\ldots,x_p\}$ with leading terms $lt(g_i)=a_ix^{\alpha_i}$ for $i=1,2$, the \emph{minimal S-process} of $g_1,g_2$ is
$$
S_{min}(g_1,g_2):=\frac{lcm(x^{\alpha_1},x^{\alpha_2})}{x^{\alpha_1}}g_1-\frac{a_1}{a_2}\frac{lcm(x^{\alpha_1},x^{\alpha_2})}{x^{\alpha_2}}g_2,
$$
where $lcm(x^{\alpha_1},x^{\alpha_2})$ is the usual least common multiple. Notice that if $\alpha_1\mid \alpha_2$ (resp. $\alpha_2\mid \alpha_1$), then we have that $S_{min}(g_1,g_2)$ is a reduction of $g_2$ modulo $\{g_1\}$ (resp. a reduction of $g_1$ modulo $\{g_2\}$).

Now, we recall an algorithm for computing a standard basis.

\bigskip 
\textbf{Büchberger's Algorithm}
\bigskip 

INPUT: $(g_1,\ldots,g_j)=I\subset \mathbb C\{x_1,x_2,\ldots,x_p\}$ ideal and a monomial order $\preceq $.

OUTPUT: $B$ standard basis of $I$. 

START: 

We put $B=\{g_1,\ldots,g_j\}$.

loop\{

\hskip 1cm for all distinct pairs of  elements $h_1,h_2\in B$:

\hskip 2cm Compute $s=S_{min}(h_1,h_2)$ and $r_\infty$ a final reduction of $s$ modulo $B$.

\hskip 2cm if $r_\infty\neq 0$:

\hskip 3cm add $r_\infty$ to $B$.

\hskip 1cm if all the final reductions computed  are 0:

\hskip 2cm Return.

\} end loop 

\bigskip 

We can use minimal standard bases for computing codimension of ideals. Consider $I\subset \mathbb C\{x_1,\ldots,x_p\}$ an ideal and the $\mathbb C$-vector space $Q=\mathbb C\{x_1,\ldots,x_p\}/I$, when the complex dimension of $Q$ is finite we have that:
$$
dim_\mathbb C Q=\#((\mathbb Z_{\geq 0})^p\setminus lp(I)).
$$

The previous formula can be computed by means of a minimal standard basis of $I$. In the particular case where $p=2$, the computation is done as follows. Assume that $B=\{g_1,\ldots,g_j\}$ is a minimal standard basis of an ideal $I\subset \mathbb C\{x,y\}$. Put $lp(g_i)=(a_i,b_i)$ for $i=1,\ldots,j$, and assume that they are ordered such that:
$$
0\leq a_1<a_2<\ldots<a_j;\quad b_1>b_2>\ldots>b_j\geq 0.
$$ 
\begin{proposition}\label{prop:pr:codimension}
With the notations as above, the dimension of $Q=\mathbb C\{x,y\}/I$ as $\mathbb C$-vector space is finite if and only if $a_1=b_j=0$. Additionally, in the finite case, we have that
$$
dim_\mathbb C Q=\sum_{i=2}^jb_{i-1}(a_i-a_{i-1}).
$$ 
\end{proposition}

The previous proposition will be used when computing intersection multiplicities.

\begin{remark}
Consider an ideal $I\subset \mathbb C\{x_1,\ldots,x_p\}$, and assume that $\{g_1,\ldots,g_j\}$ is a minimal standard basis of $I$. The set of leading powers $\{lp(g_i)\}_{i=1,\ldots,j}$ does not depend on the minimal standard basis chosen. In particular, we also have that two minimal standard basis of an ideal $I$ have the same number elements. Moreover, any ideal admits a finite standard basis.  
\end{remark}

\subsection{Plane Cusps}\label{subsec:cuspidal}
From now on, we work with the ring $\mathbb C\{x,y\}$. Let $C$ be a plane branch in $(\mathbb C^2,\mathbf 0)$ defined by an implicit equation $f=0$, or just $f$, where $f$ is a reduced irreducible element in $\mathbb C\{x,y\}$. By the Newton-Puiseux Theorem there exists a parametrization $\phi(t)=(x(t),y(t))$ of $C$, such that with $x(t),y(t)\in\mathbb C\{t\}$, see \cite{wall}. We assume that the parametrization $\phi(t)$ is primitive, that is,
$$
\nu_\mathbf 0(f)=\min\{\operatorname{ord}_t(x(t)), \operatorname{ord}_t(y(t))\}.
$$
Here $\nu_\mathbf 0(f)$ is the multiplicity of $f$ at the origin, and $\operatorname{ord}_t$ stands for the order in the variable $t$.  Any parametrization satisfying the previous condition is called \emph{primitive}. All the parametrizations appear in this work are assumed to be primitive.

Given $h_1,h_2\in\mathbb C\{x,y\}$, their \emph{intersection multiplicity} at the origin is
$$
i_\mathbf 0(h_1,h_2)=dim_\mathbb C\frac{\mathbb C\{x,y\}}{(h_1,h_2)}.
$$
The intersection multiplicity can be computed by obtaining a minimal standard basis of the ideal $(h_1,h_2)$. In the particular case where $h_2=f$, then
$$
i_\mathbf 0(f,h_1)=\operatorname{ord}_t(h_1\circ \phi).
$$
We usually denote this as  $\nu_C(h_1):=\operatorname{ord}_t(h_1\circ \phi)$. 

The \emph{semigroup} $\Gamma_C$ of the branch $C$ is defined as the set of all possible multiplicity intersections of $f$ with any other element in $\mathbb C\{x,y\}$, i.e.,
$$
\Gamma_C:=\{\nu_C(h):h\in\mathbb C\{x,y\}\},\quad \text{ where }\nu_C(h):=\operatorname{ord}_t(h\circ \phi).
$$
The semigroup satisfies that the cardinality of $\mathbb Z_{\geq 0}\setminus \Gamma_C$ is finite. Equivalently, there exists an element $c_\Gamma\in \Gamma_C$, such that $c_\Gamma$ is the minimum element in $\Gamma_C$ with the property that for any $k\geq c_\Gamma$, we have that $k\in\Gamma_C$. The element $c_\Gamma$ is called the \emph{conductor} of $\Gamma_C$.

A branch $C$ is said to be a \emph{cusp} with Puiseux pair $(n,m)$ if its semigroup $\Gamma_C$ is generated by  two integer numbers $n$ and $m$, with $2\leq n<m$ and $gcd(n,m)=1$. In other words, $\Gamma_C=\langle n,m\rangle$. For this particular case where $\Gamma_C=\langle n,m\rangle$, then the conductor of $\Gamma_C$ is 
\begin{equation}\label{eq:conductor}
c_\Gamma=(n-1)(m-1).
\end{equation}
According to \cite{wall}, for any cusp $C$, there is a system of local analytic coordinates $(x,y)$ in $(\mathbb C^2,\mathbf 0)$, such that we can find a primitive parametrization of the form
\begin{equation}\label{eq:pr:adaptadas1}
\phi(t)=(t^n,\mu_1 t^m+h.o.t.),\; \text{ with }\; \mu_1\neq 0.
\end{equation}
Equivalently, in these coordinates $(x,y)$ we can find a implicit equation $f$ of $C$ given by
\begin{equation}\label{eq:pr:adaptadas2}
f=\mu_2 x^m+y^n+\sum_{\substack{\alpha,\beta\geq 0\\ n\alpha+m\beta>nm}}z_{\alpha\beta}x^\alpha y^\beta,\; \text{ with } \mu_2\neq 0\text{ and }z_{\alpha\beta}\in\mathbb C.
\end{equation}
A system of local coordinates, such the branch $C$ satisfies one of the equivalent conditions \eqref{eq:pr:adaptadas1} or \eqref{eq:pr:adaptadas2} is said to be a local system of \emph{adapted coordinates} to $C$.
\begin{remark}\label{rem:pr:monomios-f}
In the Equation \eqref{eq:pr:adaptadas2} there are no terms of the form $xy^k$ with $k<n$.
\end{remark}

\subsection{Semimodule of Differential Values and Cuspidal Semimodules}\label{subsec:semimodule cuspidal}
All the concepts in this subsection can be found, in more detail, in \cite{yo2}. Let $C$ be a branch and $\phi(t)$ a primitive parametrization of $C$. Consider $\Omega_{\mathbb C^2,\mathbf 0}^1$ the module of holomorphic differential 1-forms in $(\mathbb C^2,\mathbf 0)$. We recall that, after choosing a system of coordinates $(x,y)$ in $(\mathbb C^2,\mathbf 0)$, the module $\Omega_{\mathbb C^2,\mathbf 0}^1$ is a free $\mathbb C\{x,y\}$-module of rank two, generated by $dx$ and $dy$. The \emph{semimodule of differential values} $\Lambda_C$ of $C$ is defined by 
$$
\Lambda_C:=\{\nu_C(\omega):\omega\in \Omega_{\mathbb C^2,\mathbf 0}^1\}.
$$
Here $\nu_C(\omega)$ is the \emph{differential value} of $\omega$ which is computed as $\nu_C(\omega):=\operatorname{ord}_t(\alpha(t))+1$, where $\phi^*(\omega)=\alpha(t)dt$ and $\phi^*\omega$ denotes the pull-back of $\omega$ by the parametrization $\phi$. 

The semimodule of differential values $\Lambda_C$ of a branch $C$ satisfies the following properties:
\begin{itemize}
\item $\Lambda_C$ is a an analytic invariant of $C$, that is to say, the semimodule of differential values does not depend on the local system of analytic coordinates chosen, see \cite{hefez2}.
\item Given a function $h$, such that $h(\mathbf 0)=0$, then $\nu_C(h)=\nu_C(dh)$. This implies that $\Gamma_C\setminus\{0\}\subset \Lambda_C$. Since the $\mathbb Z_{\geq0}\setminus \Gamma_C$ is a finite set, then $\Lambda_C\setminus \Gamma_C$ is also finite. Moreover, every element $\lambda\in (\Lambda_C\setminus \Gamma_C)$ satisfies that $\lambda<c_\Gamma$, where $c_\Gamma$ is the conductor of $\Gamma_C$.
\item For any $h\in\mathbb C\{x,y\}$ and $\omega\in\Omega_{\mathbb C^2,\mathbf 0}^1$, then $\nu_C(h\omega)=\nu_C(h)+\nu_C(\omega)$ . This implies that, for any $\gamma\in\Gamma_C$ and $\lambda\in \Lambda_C$, then we have that $\lambda+\gamma\in\Lambda_C$.
\end{itemize}  
The third property is equivalent to say that $\Lambda_C$ is a $\Gamma_C$-semimodule. More generally, consider $\Gamma\subset\mathbb Z_{\geq0}$ a monoid. We say that set $\Lambda\subset \mathbb Z_{\geq0}$ is a $\Gamma$-semimodule if the following property holds: for all $\gamma\in\Gamma$ and $\lambda\in\Lambda$, we have that $\lambda+\gamma\in\Lambda$. It comes through this notion that we can study some of the underlying structure of the semimodule of differential values $\Lambda_C$.

In general, for any $\Gamma$-semimodule $\Lambda$, we have that there exists a unique increasing sequence $(\lambda_{-1},\lambda_0,\lambda_1,\ldots,\lambda_s)$ of elements in $\Lambda$, where $s\geq -1$ is minimal satisfying the property:
$$
\Lambda=\bigcup_{i=-1}^s(\lambda_i+\Gamma).
$$
The sequence $(\lambda_{-1},\lambda_0,\ldots,\lambda_s)$ is called the \emph{basis} of the semimodule $\Lambda$, each element $\lambda_i$ of the basis is called the $i-$\emph{element of the basis}; the parameter $s$ is said to be the length of the basis. By definition of a basis, we have that $\lambda_i\notin \lambda_j+\Gamma$ for $i\neq j$. Therefore, if we denote by $n$ the minimum element in $\Gamma$, $n=\min (\Gamma)$, then we have that $\lambda_i$ and $\lambda_j$ define different classes modulo $n$, this implies that $s\leq n-2$. In the case of a semimodule of differential values of a singular branch $C$, then $s\geq 0$, where $s=0$ if and only $C$ is quasihomogeneous, see \cite{hefez2}. 

Consider $(\lambda_{-1},\lambda_0,\ldots,\lambda_s)$ the basis of $\Lambda_C$, by definition of the semimodule of differential values $\Lambda_C$, there exists a sequence of 1-forms $(\omega_{-1},\omega_0,\ldots,\omega_s)$ such that $\nu_C(\omega_i)=\lambda_i$ for $i=-1,0,\ldots,s$. The sequence $(\omega_{-1},\ldots,\omega_s)$ is called a \emph{minimal standard basis of the module of differentials} of $C$. When there is no confusion, we just call it minimal standard basis of $C$.
\begin{remark}
Saying $(\omega_{-1},\omega_0,\ldots,\omega_s)$ is a minimal standard basis of the module of differentials is an abuse of notation. In fact, $(\phi^*(\omega_{-1}),\phi^*(\omega_0),\ldots,\phi^*(\omega_s))$ is truly a minimal standard basis of the module of differentials of $C$ in the sense of \cite{hefez1}. The second definition corresponds to a closer notion of standard basis to the one given in section \ref{sec:ideal}, but for the case of modules instead of ideals. 
\end{remark}

Given $(\lambda_{-1},\ldots,\lambda_s)$ the basis of a $\Gamma$-semimodule $\Lambda$. We define the \emph{decomposition sequence} of $\Lambda$ as:
$$
\Lambda_{-1}\subset \Lambda_0\subset \ldots\subset \Lambda_s=\Lambda\;\text{ with }\; \Lambda_{i}=\bigcup_{j=-1}^i(\lambda_j+\Gamma),\; i=-1,0,\ldots,s.
$$
\begin{remark}\label{rem:pr:base}
We have that $\lambda_{-1}=\min(\Lambda)$ and $\lambda_i=\min (\Lambda\setminus \Lambda_{i-1})$ for $i=0,1,\ldots,s$.
\end{remark}

A $\Gamma$-semimodule $\Lambda$ is said to be a \emph{cuspidal semimodule} if $\Gamma$ is generated by a pair $(n,m)$, with $2\leq n<m$ and $gcd(n,m)=1$, that is, the underlying semigroup corresponds with one associated to a cusp.

If $\Lambda$ is a cuspidal $\Gamma$-semimodule with $\Gamma=\langle n,m\rangle$ and basis $(\lambda_{-1},\lambda_0,\ldots,\lambda_s)$, we define the \emph{axes} $u_i$ for $i=1,2,\ldots,s+1$ as:
$$
u_i=\min((\lambda_{i-1}+\Gamma)\cap \Lambda_{i-2}).
$$
We say that $\Lambda$ is an \emph{increasing semimodule} if $\lambda_i>u_i$ for $i=1,2,\ldots,s$, see \cite{alberich}. 
\begin{remark}\label{rem:pr:escritura ui}
In \cite{yo} it was shown that for each axis $u_i$ with $i=1,\ldots,s$, there exists a unique $k<i-1$, such that 
$$
u_i=\lambda_{i-1}+\gamma_1=\lambda_k+\gamma_2\;\text{ with }\; \gamma_1,\gamma_2\in \Gamma\setminus\{0\}.
$$
Besides, if we write $\gamma_1=na_1+mb_1$ and $\gamma_2=na_2+mb_2$, it is satisfied that either $b_1=a_2=0$ or $a_1=b_2=0$.
\end{remark}
From the axes we can define the \emph{critical values} $t_i$ for $i=-1,0,\ldots,s+1$ as follows:
\begin{itemize}
\item First, we put $t_{-1}=\lambda_{-1}$ and $t_0=\lambda_0$.
\item Next, for $j=1,2,\ldots,s+1$, we define inductively $t_j=t_{j-1}+u_j+\lambda_{j-1}$. 
\end{itemize}

\subsection{Semimodule of differential values for cusps}\label{subsec:cusps}

Fix $C$ a cusp with Puiseux pair $(n,m)$, semigroup $\Gamma_C=\langle n,m\rangle$ and with semimodule of differential values $\Lambda_C$. By definition $\Lambda_C$ is a cuspidal semimodule. Denote by $(\lambda_{-1},\lambda_0,\ldots,\lambda_s)$ the basis of $\Lambda_C$ and consider the axes and critical values $(u_1,u_2,\ldots,u_s)$ and $(t_{-1},t_0,\ldots,t_s)$. In \cite{delorme}, it is shown that the semimodule of differential values  $\Lambda_C$ satisfies the following conditions:
\begin{itemize}
\item $\lambda_{-1}=n$ and $\lambda_0=m$. 
\item $\Lambda_C$ is increasing semimodule, equivalentely, that $\lambda_i>u_i$ for $i=1,2,\ldots,s$.
\end{itemize}
Moreover, in \cite{alberich} it is proven the inverse. Take  $\Lambda$ a cuspidal $\Gamma_C$-semimodule. Then $\Lambda$ is the semimodule of differential values of a cusp $D$ if and only if: $\Lambda$ is increasing and the (-1)-element and the 0-element of the basis are $n$ and $m$ respectively. 

We want to enphasize that $\lambda_1-n$ is the well known Zariski's invariant, see \cite{zariski}.

\begin{remark}\label{rem:pr:u1 y t1}
Since $\lambda_{-1}=n$ and $\lambda_0=m$, we always have that both the axis $u_1$ and the critical value $t_1$ are equal to $n+m$. Moreover, when considering the decomposition sequence $\Lambda_{-1}\subset \Lambda_0\subset\ldots\subset \Lambda_s=\Lambda_C$ of $\Lambda_C$, we have that $\Gamma_C=\Lambda_0\cup \{0\}$.
\end{remark}

From the results in \cite{yo,delorme}, it was obtained a method to compute a minimal standard basis of the module of differentials for a cusp. We proceed to detail the algorithm. But first, we introduce the following notation that we use along the rest of the paper.
\begin{notation}
Given two 1-forms $\omega,\omega'\in\Omega_{\mathbb C^2,\mathbf 0}^1$, such that $\nu_C(\omega)=\nu_C(\omega')<\infty$, we denote by $\mu^+$ the unique non zero constant, such that $\nu_C(\omega+\mu^+\omega')>\nu_C(\omega)$. This constant $\mu^+$ is called the \emph{tuning constant}.
\end{notation}
Consider $f_{-1},f_0\in\mathbb C\{x,y\}$ such that $\nu_C(f_{-1})=n$ and $\nu_C(f_0)=m$. Take $\omega,\omega'\in\Omega_{\mathbb C^2,\mathbf 0}^1$ with $\nu_C(\omega)\neq \infty \neq \nu_C(\omega')$. Then $\omega'$ is a said to be a \emph{reduction} of $\omega$ modulo $H\subset \Omega_{\mathbb C^2,\mathbf 0}^1$, if there exist $\eta\in H$, $a,b\in\mathbb Z_{\geq 0}$, such that 
$$
\nu_C(\omega)=f_{-1}^af_0^b\eta\text{ with }\omega'=\omega+\mu^+f_{-1}^af_0^b\eta.
$$ 
In other words, we have that $\nu_C(\omega')>\nu_C(\omega)$. 

Furthermore, we can consider a sequence of reductions starting from $\omega$, and in the limit, we obtain an element $r_\infty$ satisfying that either $\phi^*r_\infty$ is zero or is no longer reducible modulo $B$. This element $r_\infty$ is called a \emph{final reduction} of $\omega$ modulo $B$.

\bigskip 

\textbf{Delorme's Algorithm}
\bigskip

INPUT: A cusp $C$ and $f_{-1},f_0$ with $\nu_C(f_{-1})=n$ and $\nu_C(f_0)=m$. 

OUTPUT: $H$ minimal standard basis of the module of differentials of $C$.

START: 

Put $H=\{\omega_{-1}=df_{-1},\omega_0=df_0\}$, $i=1$,  $B_0=(\lambda_{-1}=n,\lambda_0=m)$ and $\Lambda_{-1}=(n+\Gamma)$.

loop\{

\hskip 1cm Compute $u_i=\min\{(\lambda_{i-1}+\Gamma)\cap \Lambda_{i-2}\}$ and put $\Lambda_{i-1}=\Lambda_{i-2}\cup (\lambda_{i-1}+\Gamma)$.

\hskip 1cm Put $\eta_1=h\omega_{i-1}$ and $\eta_2=h'\omega_k$ with $k<i-1$. Here $h=f_{-1}^af_0^b$, $h'=f_{-1}^{a'}f_0^{b'}$ with $a,a',b,b'\in\mathbb Z_{\geq 0}$ and such that $\nu_C(h\omega_{i-1})=\nu_C(h'\omega_k)=u_i$. 

\hskip 1cm Put $\eta=\eta_1+\mu^+\eta_2$.

\hskip 1cm Compute $r_\infty$ a final reduction of $\eta$ modulo $H$.

\hskip 1cm if $\nu_C(r_\infty)=\infty$:

\hskip 2cm Return.

\hskip 1cm otherwise:

\hskip 2cm Add $\omega_i=r_\infty$ to $H$.

\hskip 2cm Put $\lambda_{i}=\nu_C(\omega_i)$, $B_i=(n,m,\lambda_1,\ldots,\lambda_i)$ and $i=i+1$.

\} end loop 

\bigskip

The previous algorithm always returns a minimal standard basis. Moreover, in each iteration of the loop, it computes a single 1-form. In contrast, the algorithm presented in \cite{hefez1} offers a method to compute standard bases for broader families of curves. However, it may compute additional 1-forms. In \cite{yo2}, it was shown how Delorme's algorithm can be adapted to compute a Saito basis of $C$, that is, a basis of the module of 1-forms with $C$ invariant.

In an adapted system of coordinates $(x,y)$, the 1-forms of a minimal standard basis of a cusp can be characterised in term of a monomial value. Take the pair $(n,m)$ and assume that the leading power of an element $g\in\mathbb C\{x,y\}$ is $lp(g)=(a,b)$. We define the \emph{monomial value} of $g$ with weights $(n,m)$ as 
$$
\nu_{n,m}(g):=na+mb.
$$
We extend the notion of monomial value to 1-forms as follows: given $\omega\in\Omega_{\mathbb C^2,\mathbf 0}^1$ with $\omega=Adx+Bdy$, the monomial value of $\omega$ is 
$$\nu_{n,m}(\omega):=\min\{\nu_{n,m}(xA),\nu_{n,m}(yB)\}.
$$ 
In fact, as shown in \cite{yo}, the previous monomial value can be interpreted in terms of divisorial valuations. These divisorial valuations make the results independent of the coordinate system. However, for the purposes of this work, we do not need such a generality.

\begin{remark}\label{rem:pr:xy}
In adapted coordinates $(x,y)$, we can assume that $f_{-1}=x$ and $f_0=y$. In general, we have that $\{\omega_{-1},\omega_0\}$ is a basis of $\Omega_{\mathbb C^2,\mathbf 0}^1$, see \cite{yo}. Besides, their monomial values are fixed. In fact, we have that $\nu_{n,m}(\omega_{-1})=\lambda_{-1}=t_{-1}=n$ and $\nu_{n,m}(\omega_0)=\lambda_0=t_0=m$. More precisely, in adapted coordinates the 1-forms can be written as $\omega_{-1}=\mu dx+\omega$ and $\omega_0=\mu' dy+\omega'$, where $\mu,\mu'$ are non zero constants and; the monomial values of $\omega$ and $\omega'$ satisfy that $\nu_{n,m}(\omega)>\nu_{n,m}(dx)=n$ and $\nu_{n,m}(\omega')>\nu_{n,m}(dy)=m$.
\end{remark}

The previous remark describes the first two 1-forms of a minimal standard basis of a cusp. For the other ones, we have the following characterization.

\begin{theorem}[\cite{yo} Theorem 7.13]\label{thm:standard:delorme}
Assume that we are in local system of adapted coordinates to a cusp $C$ with pair Puiseux $(n,m)$. Consider $(\lambda_{-1},\lambda_0,\ldots,\lambda_s)$ the basis of the semimoduple differential values $\Lambda_C$ of $C$.  For each $1\leq i\leq s$ we have the following statements:
\begin{enumerate}
\item $\lambda_i=\sup\{\nu_C(\omega):\nu_{n,m}(\omega)=t_i\}$.
\item If $\nu_C(\omega)=\lambda_i$, then $\nu_{n,m}(\omega)=t_i$.
\item For each $1$-form $\omega$ with $\nu_C(\omega)\notin \Lambda_{i-1}$, there is a unique pair $a,b\geq 0$ such that $\nu_{n,m}(\omega)=\nu_{n,m}(x^ay^b\omega_i)$. Moreover, we have that $\nu_C(\omega)\geq \lambda_i+na+mb$.
\item Let $k=\lambda_i+na+mb$, then $k\notin \Lambda_{i-1}$ if and only if for all $\omega$ such that $\nu_C(\omega)=k$ we have that $\nu_{n,m}(\omega)\leq \nu_{n,m}(x^ay^b\omega_i)$.
\item We have that $\lambda_i>u_i$.
\end{enumerate}
\end{theorem}
We recall that $u_i$ and $t_i$ are the axes and the critical values as defined in the section \ref{subsec:semimodule cuspidal}. 

\section{Semimodule of differential values from an implicit equation}\label{sec:impli}

This section is devoted to explain how to apply Delorme's algorithm without using a primitive parametrization. We approach two questions: first, how to compute a differential value using an implicit equation. Second, how to compute the tuning constants $\mu^+$ needed in several steps of the algorithm. The first question has a well known solution. Take $C$ a branch defined by the implicit equation $f=0$. Consider $\omega=Adx+Bdy$ a 1-form, we denote by $X_\omega\in\mathcal X_{\mathbb C^2,\mathbf 0}$ to the vector field 
$$
X_\omega:=B\partial_x-A\partial_y.
$$ 
Notice that the definition of $X_\omega$ depends on the chosen coordinate system. We have the following result:
\begin{lemma}[\cite{corral} Proposition 2.3]\label{lem:impli:valor-diferencial}
Let $C$ be a branch (not necessarily a cusp) defined by the implicit equation $f=0$. Then we have that for any 1-form $\omega\in\Omega_{\mathbb C^2,\mathbf0}^1$, the differential value $\nu_C(\omega)$ is given by
\begin{equation}\label{eq:impli:valor-diferencial}
\nu_C(\omega)=i_\mathbf 0(X_\omega(f),f)-c_\Gamma+1.
\end{equation}
\end{lemma}
In fact, the previous lemma is a weaker version of the Proposition 2.3 in \cite{corral}.

Notice that $i_\mathbf 0(X_\omega(f),f)$ can be computed using a minimal standard basis of the ideal $(X_\omega(f),f)$, which is obtained by using an implicit equation. Hence, we can find the differential value $\nu_C(\omega)$ without directly using a parametrization.

Now assume that $C$ is a cusp with Puiseux pair $(n,m)$ and $(x,y)$ is a local system of coordinates adapted to $C$. Take an implicit equation $f=0$ of $C$ as in Equation \eqref{eq:pr:adaptadas2}, that is
$$
f=\mu x^m+y^n+\sum_{\substack{\alpha,\beta\geq 0\\ n\alpha+m\beta>nm}}z_{\alpha\beta}x^\alpha y^\beta,\; \text{ with } \mu\neq 0\text{ and }z_{\alpha\beta}\in\mathbb C.
$$
With this setting, finding the tuning constants $\mu^+$ relies on computing a minimal standard basis of the ideal $(X_\omega(f),f)$. We want to remark that the procedure for computing the tuning constants only works in the cuspidal case. As in Example \ref{ex:pr:orders}, we consider the weighted order with respect to $(n,m)$, where we recall that $(a,b)\prec(c,d)$ if and only if either $na+mb<nc+md$ or $na+mb=nc+md$ and $a<c$. 


\begin{proposition}\label{prop:impli:base-ideal}
Consider a vector field $X\in \mathcal X_{\mathbb C^2,\mathbf 0}$ and the ideal $I=(X(f),f)$, with the assumption that $f\nmid X(f)$.
\begin{itemize}
\item[1.] Suppose that $lp(X(f))\neq (0,n-1)$. Denote by $h$ a final reduction of $X(f)$ modulo $\{f\}$. Then:
\begin{enumerate}
\item[\textsl{a})] If $lp(h)=(a,0)$, then $\{f,h\}$ is a minimal standard basis of $I$.
\item[\textsl{b})] If $lp(h)=(a,b)$ with $b>0$, then $\{f,h,S_{min}(f,h)\}$ is a minimal standard basis of $I$, where $S_{min}(f,h)$ is a minimal $S$-process of $f$ and $h$.
\end{enumerate}
\item[2.] Suppose that $lp(X(f))=(0,n-1)$, and let $g$ be a final reduction of $f$ modulo $\{X(f)\}$. Then $\{X(f),g\}$ is a minimal standard basis of $I$, where $lp(g)=(m,0)$.
\end{itemize}
\end{proposition}

Notice that in the case \textit{1}.\textsl{b}), we are saying that $S_{min}(f,h)$ is not reducible modulo $\{f,h\}$. 
\begin{proof}
The proofs of the three statements are pretty similar; for this reason, we omit the one for Statement \textit{2}. Assume, as in Statement \textit{1}, that $lp(X(f))\neq (0,n-1)$. This condition is equivalent to $lp(X(f))\neq lp(f_y)$, that is, meaning that the term $\alpha \partial _y$ with $\alpha\in\mathbb{C}$ of the vector field $X$ is zero.  From Remark \ref{rem:pr:monomios-f}, we also deduce that $lp(X(f))\neq (0,c)$, for any $c<n$.

To compute a minimal standard basis of the ideal $I$, we apply Büchberger's algorithm. First, we find a generator system $B=\{f_1,f_2\}$ of $I$, such that neither $f_1$ is reducible by $f_2$, nor vice versa. Next, we will compute the minimal $S$-process $S_{min}(f_1,f_2)$ and a final reduction modulo $\{f_1,f_2\}$. Finally, we iterate the process as many times as needed. In fact, we show that  either $\{f_1,f_2\}$ is the desired standard basis, or the algorithm finds one after the first iteration. 

Since $f$ does not divide $X(f)$, we know that $I\neq (f)$, meaning that we need at least two generators for the minimal standard basis. Because $lp(X(f))\neq (0,c)$ with $c<n$ and $lp(f)=(0,n)$, if $f$ is reducible modulo $\{X(f)\}$, then $lp(X(f))=(0,n)$, implying that $X(f)$ is reducible modulo $\{f\}$. Therefore, we set $f_1=f$ and $f_2=h$, where $h$ is a final reduction of $X(f)$ modulo $\{f\}$.

Let the leading power of $h$ be $lp(h)=(a,b)$ with $b<n$. Otherwise, $h-\mu' x^ay^{n-b}f$ would give us a new reduction of $h$ for an appropriate constant $\mu'$, contradicting the fact that $h$ is a final reduction of $X(f)$ modulo $\{f\}$.

We have to check that for the case $b=0$, the algorithm has found a standard basis, and hence minimal, since neither $f$ nor $h$ can be omitted. On the contrary, if $b\neq 0$, we have to compute an extra element. In both cases, we can write
$$
S_1=S_{min}(f,h)=x^af-\mu_1 y^{b-n} h,
$$
where $\mu_1$ is the unique constant such that $lp(S_1)\succ lp(x^af)=lp(y^{n-b}f)$.

We claim that the leading term of $S_1$ is $lt(S_1)=\mu x^{a+m}$, where $\mu$ is the constant appearing in the implicit equation $f$ as in Equation \eqref{eq:pr:adaptadas2}. This is because of the following two facts: first, the Newton cloud of $h$ has no points of the form $(c,d)$ with $nc+md=na+mb$ and $d\geq n$. Otherwise, it would follow that $c<a$ and that $lp(h)\neq (a,b)$. Second, the leading term of $x^af-lt(x^af)$ is $\mu x^{a+m}$. Since $lt(S_1)=\mu x^{a+m}$, then $S_1$ does not admit any reduction modulo $\{f,h\}$, unless $b=0$.

\emph{Case} $b=0$: for any partial reduction $r$ of $S_1$, we have that 
$$
lp(r)=(c,d)\succ lp(S_1)=(a+m,0)\succ (a,n)=lp(x^af).
$$ 
We have the following:
\begin{itemize}
\item If $d\geq n$, then $r$ can be reduced modulo $\{f\}$. 
\item If $c\geq a$, then it can be reduced modulo $\{h\}$.
\item If  $c<a$ and $d<n$, then  we have that $nc+md<na+mn$, in contradiction with the assumption that $(c,d)\succ (a,n)$. 
\end{itemize}
In conclusion, if $b=0$, then $r$ is reducible modulo $\{f,h\}$.  Since $r$ is any partial reduction of $S_1$ modulo $\{f,h\}$, it follows that $0$ is a final reduction of $S_1$ modulo $\{f,h\}$. Thus the Büchberger's algorithm stops and $\{f,h\}$ is a minimal standard basis of $I$. This ends the proof of Statement \textit{1}.\textsl{a}).

\bigskip

\emph{Case} $b\neq 0$: As we said before, since $lt(S_1)=\mu x^{a+m}$, we have that in this case $S_1$ is its own final reduction modulo $\{f,h\}$. Therefore, by the Büchberger's algorithm, it is necessary to add $S_1$ to our candidate of standard basis $\{f,h\}$. We only need to verify that the algorithm stops here. In other words, we have to compute all new possible minimal $S$-process and see that they have 0 as a final reduction modulo $\{f,h,S_1\}$. 

There are only two new minimal $S$-process to consider:
\begin{eqnarray*}
S_2=S_{min}(f,S_1)&=&x^{a+m}f-\mu_2 y^nS_1\\
S_3=S_{min}(h,S_1)&=&x^{m}h-\mu_3 y^{b}S_1.
\end{eqnarray*} 

Here $\mu_2$ and $\mu_3$ are the unique constants such that $lp(S_2)\succ lp(x^{a+m}f)=(a+m,n)$ and $lp(S_3)\succ (x^mh)=(a+m,b)$. In order to show that their final reductions modulo $\{f,h,S_1\}$ are zero, it suffices to show that any pair $(c,d)\succ (a+m,b)$ is divisible by $(0,n),(a,b)$ or $(a+m,0)$, that is, the leading powers of $f,h$ and $S_1$.

We have that $nc+md\geq n(a+m)+mb$, hence 
\begin{itemize}
\item If $0\leq c\leq a$, then $d\geq n+b>n$ and $(0,n)$ divides $(c,d)$. 
\item If $a<c<m+a$, then $d>b$ and $(a,b)$ divides $(c,d)$.
\item If $c\geq a+m$, then $(a+m,0)$ divides $(c,d)$. 
\end{itemize}
Hence the final reductions of $S_2$ and $S_3$ modulo $\{f,h,S_1\}$ are 0, ending the proof of Statement \textit{1}.\textsl{b}).
\end{proof}



\begin{proposition}\label{prop:impli:multiplicidad}
Take a 1-form $\omega\in\Omega_{\mathbb C^2,\mathbf 0}^1$, such that $f\nmid X_\omega(f)$. Put $h\in\mathbb C\{x,y\}$ a final reduction of $X_\omega(f)$ modulo $\{f\}$, with leading power $lp(h)=(a,b)$.  Then we have that
$$
\nu_C(\omega)=n(a+1)+m(b+1)-nm.
$$
\end{proposition}
\begin{proof}
By Lemma \ref{lem:impli:valor-diferencial} we have that
$$
\nu_C(\omega)=i_{\mathbf 0}(X_\omega(f),f)-c_\Gamma+1.
$$
By Propositions \ref{prop:pr:codimension} and \ref{prop:impli:base-ideal}, we find that $i_{\mathbf 0}(X_\omega(f),f)=na+mb$. We conclude by noting that $c_\Gamma=(n-1)(m-1)$.
\end{proof}

Proposition \ref{prop:impli:multiplicidad} solves the problem of computing the tuning constants $\mu^+$, at least in the case we are interested in. Suppose that $\eta_1,\eta_2\in\Omega_{\mathbb C^2,\mathbf 0}^1$ are two 1-forms, with $\nu_C(\eta_1)=\nu_C(\eta_2)<\infty$. Let us see how to compute the tuning constant of the 1-form $\eta_1+\mu^+\eta_2$. First, consider $h_1$ and $h_2$ to be final reductions of $X_{\eta_1}(f)$ and $X_{\eta_2}(f)$ modulo $\{f\}$ respectively. Put $lp(h_1)=(a,b)$ and $lp(h_2)=(c,d)$ the leading powers, where $b,d<n$. By Proposition \ref{prop:impli:multiplicidad}, we have that 
$$
n(a+1)+m(b+1)-nm=\nu_C(\eta_1)=\nu_C(\eta_2)=n(c+1)+m(d+1)-nm.
$$
Because of the inequalities $0\leq b,d<n$, we find that $(a,b)=(c,d)$. The previous equalities can be translated in terms of leading terms as
$$
lt(h_1)=\mu_1 x^ay^b;\quad lt(h_2)=\mu_2 x^ay^b.
$$
Here $\mu_1$ and $\mu_2$ are non zero constants. If we put $\mu^+=-\mu_1/\mu_2$ as the tuning constant, we observe the following: the leading power $lp(h_1+\mu^+h_2)\succ lp(h_1),lp(h_2)$. Since  
$h_1+\mu^+h_2$ is, at least, a partial reduction of $X_{\eta_1+\mu^+\eta_2}(f)$, by Propositions \ref{prop:impli:multiplicidad}, we verify that $\nu_C(\eta_1+\mu^+\eta_2)>\nu_C(\eta_1),\nu_C(\eta)$. 


Since we know how to compute differential values and tuning constants, we can apply the Delorme's algorithm without a parametrization.


\section{Bernstein-Sato Polynomial}

Consider the ring of non-commutative power series $A=\mathbb{C}\{x_1,\ldots,x_p,\partial_1,\ldots,\partial_p\}$ in $2p$ variables, and define $\mathcal{D}$ to be the quotient of $A$ by the commutators $[x_i,x_j]=0$ and $[\partial_i,x_j]=\delta_{ij}$, where $\delta_{ij}$ is the Kronecker delta. The ring $\mathcal{D}$ is the set of differential operators in $p$ variables, whose action on $\mathbb{C}\{x_1,\ldots,x_p\}$ is defined in the natural way.

We take $\mathcal D[\rho]$ the ring of polynomials in the variable $\rho$ and coefficients in $\mathcal{D}$. For any function $g\in \mathbb{C}\{x_1,\ldots,x_p\}$, we can extend the action of $\mathcal{D}$ to functions of the form $g^\rho$ by defining $\partial_i g^\rho=\rho g^{\rho-1}\partial_i g$.

Take an hypersurface $H$ of $(\mathbb C^p,\mathbf 0)$ defined by an implicit equation $g=0$ with $g\in\mathbb C\{x_1,\ldots,x_p\}$. According to \cite{bjork}, there exist non zero elements $P\in \mathcal D[\rho]$ and $B(\rho)\in\mathbb C[\rho]$ such that
$$
Pg^{\rho+1}=B(\rho)g^\rho.
$$
The set of all $B(\rho)$ for which there is an operator $P\in \mathcal D[\rho]$ satisfying the last condition is an ideal in $\mathbb C[\rho]$. Hence, it is principal, the monic generator is denoted by $b(\rho)$ and is called the Bernstein-Sato polynomial of $g$ or $H$. The Bernstein-Sato polynomial is an analytic invariant of $H$, see \cite{artal}.

Now we particularize to the case of cusps. For the rest of this section, we fix $C$ to be a cusp with Puiseux pair $(n,m)$. Many results of this work require the use of an adapted system of coordinates with respect to $C$. Nonetheless, when dealing with the Bernstein-Sato polynomial we to use more particular coordinate systems. We refer to them as nice coordinates, and we proceed to define them.

To define nice coordinates, we first introduce what we call the \emph{cuspidal sets} $P,J,M$. These sets usually appear when studying cusps, see \cite{cassou,peraire} for more details. They are defined as follows:
\begin{eqnarray*}
P&:=&\{(p_{1},p_{2})\in (\mathbb Z_{\geq 0})^2:  0\leq p_1<m-1,0\leq p_2<n-1,\; np_1+mp_2>nm\},\\
J&:=&\{ j=p_{1,j}n+p_{2,j}m-nm: (p_{1,j},p_{2,j})\in P\},\\
M&:=&\{(m-p_1-1,n-p_2-1):(p_1,p_2)\in P\}.
\end{eqnarray*}
Given $j\in J$ when we write $(p_{1,j},p_{2,j})$, we are referring to the unique element in $P$ such that $j=p_{1,j}n+p_{2,j}m-nm$. Note that if $n=2$ the three sets are empty.
\begin{remark}\label{rem:pr:na+mb=c}
Given $na+mb=na'+mb'=c$, with $a,a',b,b',c$ natural numbers, if $c<nm$, then we have that $(a,b)=(a',b')$.
\end{remark}
There exists a system of \emph{nice coordinates} such that $C$ has an implicit equation $f$ given by
\begin{equation}\label{eq:pr:nice-equation}
f=x^m+y^n+\sum_{j\in J} z_j x^{p_{1,j}}y^{p_{2,j}};\quad z_j\in \mathbb C,
\end{equation}
(see \cite{cassou}). An implicit equation like the one of the Equation \eqref{eq:pr:nice-equation} will be called a \emph{nice equation} of $C$. We fix $f\in\mathbb{C}\{x,y\}$ a nice equation of $C$ written as in Equation \eqref{eq:4:nice}.

From the cuspidal set $J$ we can define candidates for being roots of the Bernstein-Sato of $C$. Take $j\in J$, we define the rational numbers 
$$
\alpha_j=\frac{np_{1,j}+mp_{2,j}+n+m}{nm}; \quad \beta_j=\alpha_j-1=\frac{j+n+m}{nm}.
$$
We know that either $-\alpha_j$ or $-\beta_j$ is a root of the Bernstein polynomial of $C$.

According to \cite{cassou}, for any $(a,b)\in M$, it can be defined a complex function $I_0((a,b),f)(\rho)$, such that its residue at the point $\rho=-\beta_j$ is $Res_f(a,b)(\beta_j)$ and satisfies the following equation.
\begin{equation}\label{eq:residuos}
Res_f(a,b)(\beta_j)=\frac{\Gamma(\beta_j)^{-1}}{nm}\sum_{\substack{(\delta_\ell)_{\ell\in J}\in (\mathbb Z_{\geq 0})^{\# J}\\ \sum\delta_\ell\ell=k}} (-1)^{\sum \delta_\ell}\Gamma\Big( \tfrac{\sum\delta_\ell p_{1,\ell}+a}{m}\Big)\Gamma\Big( \tfrac{\sum\delta_\ell p_{2,\ell}+b}{n}\Big)\prod \frac{z_\ell^{\delta_\ell}}{\delta_\ell!},
\end{equation}
where $k=\beta_jnm-na-mb$ and $\Gamma(-)$ is the Euler's Gamma function. Also in \cite{cassou} the author shows the following result:
\begin{theorem}\label{thm:roots}
Assume that $C$ is defined by the nice equation
$$
f=x^m+y^n+\sum_{j\in J}z_j x^{p_{1,j}}y^{p_{2,j}}.
$$
and consider $b(\rho)$ its Bernstein-Sato polynomial. Then $-\beta_j$ is a root of $b(\rho)$ if and only if there exists $(a,b)\in M$, such that $Res_f(a,b)(\beta_j)\neq 0$. Otherwise, if $-\beta_j$ is not a root of $b(\rho)$, then $-\alpha_j$ is.
\end{theorem}

As mentioned in the introduction, the main goal of this work is to proof the following two results.

\begin{theorem*}[\ref{thm:bern:l1}]
Let $C$ be a cusp with semigroup $\Gamma_C=\langle n,m\rangle$ and semimodule of differential values $\Lambda_C$. Assume that $\lambda_1=\min(\Lambda_C\setminus\Gamma_C)$ exists. Then for any element $\lambda\in (\lambda_1+\Gamma_C)\setminus \Gamma_C\subset \Lambda_C$, the rational number $-\lambda/nm$ is a root of the Bernstein-Sato polynomial of $C$.
\end{theorem*}

\begin{theorem*}[\ref{thm:bern:4}]
Let $C$ be a cusp with semigroup $\Gamma_C=\langle n,m\rangle$ and semimodule of differential values $\Lambda_C$. Assume that $n\leq 4$. Then for any element $\lambda\in \Lambda_C\setminus \Gamma_C$, the rational number $-\lambda/nm$ is a root of the Bernstein-Sato polynomial of $C$.
\end{theorem*}

Before giving the proofs of both theorems, we state two useful lemmas that relate the cuspidal sets and the semimodule of differential values.

\begin{lemma}\label{lem:pr:valores no semi}
We have that
$$ 
J= \{\ell\in \mathbb{N}: \ell+n,\ell+m\notin\Gamma_C\}.
$$
As a consequence, for any $\ell\in J$, the element $\ell+n+m$ does not belong to the semigroup $\Gamma_C$.
\end{lemma}

\begin{proof}
For a proof of $J=\{\ell\in \mathbb{N}: \ell+n,\ell+m\notin\Gamma_C\},$ see \cite{peraire} Lemma 1.4.

Now, if $\ell+n+m\in\Gamma_C$, this implies that $\ell+n+m=na+mb$, with at least one of the coefficients different from 0. Hence, either $\ell+n\in \Gamma_C$ or $\ell+m\in\Gamma_C$.
\end{proof}
\begin{lemma}\label{lem:pr:valores no semi 2}
Consider a semigroup $\Gamma_C=\langle n,m\rangle$ with $3\leq n<m$ and $gcd(n,m)=1$. The following statements hold:
\begin{itemize}
\item[i)] if $\lambda\notin\Gamma_C$ with $\lambda>n+m$, then there exists $j\in J$ such that $\lambda=j+n+m$.
\item[ii)] If $\lambda+na+mb\notin\Gamma_C$ with $a,b\geq 0$, we have that the element $(a+1,b+1)\in M$.
\end{itemize}
\end{lemma} 
\begin{proof}
The first part is a consequence of Lemma \ref{lem:pr:valores no semi}, because we have that $\lambda-n-m,\lambda-n,\lambda-m\notin \Gamma_C$.

In order to prove $ii)$, consider $\lambda+na+mb\notin\Gamma_C$. Then $\lambda\notin \Gamma_C$ and by $i)$, we can write 
\begin{eqnarray*}
\lambda=j+n+m& \text{with} &j=np_{1,j}+mp_{2,j}-nm. \\
\lambda+na+mb=j'+n+m& \text{with} &j'=np_{1,j'}+mp_{2,j'}-nm.
\end{eqnarray*}
Here $(p_{1,j},p_{2,j}),(p_{1,j'},p_{2,j'})\in P$. Note that $\lambda+na+mb$ is smaller than the conductor $c_\Gamma$, that is, we have that $\lambda+na+mb<c_\Gamma=nm-n-m+1$. Since $\lambda>n+m$, we get that 
\begin{equation}\label{eq:lema1}
na+mb<nm-2n-2m.
\end{equation}
Thus, we obtain that
\begin{itemize}
\item $n(a+1)+m(b+1)<nm$.
\item $0\leq a \leq m-3;\quad 0\leq b\leq n-3$.
\end{itemize}
By Remark \ref{rem:pr:na+mb=c}, there is unique expression of $n(a+1)+m(b+1)$ given by
\begin{eqnarray*}
a+1&=&p_{1,j'}-p_{1,j}+1=m-(m-p_{1,j'}+p_{1,j}-2)-1\\
b+1&=&p_{2,j'}-p_{2,j}+1=n-(n-p_{2,j'}+p_{2,j}-2)-1.
\end{eqnarray*}
We need to check that $(m-p_{1,j'}+p_{1,j}-2,n-p_{2,j'}+p_{2,j}-2)\in P$. By definition of $P$, we have to show that:
\begin{enumerate}
\item[a)] 0 $\leq m-p_{1,j'}+p_{1,j}-2 \leq m-2$.
\item[b)] 0 $\leq n-p_{2,j'}+p_{2,j}-2 \leq n-2$.
\item[c)] $ n(m-p_{1,j'}+p_{1,j}-2)+m(n-p_{2,j'}+p_{2,j}-2)>nm$.
\end{enumerate}
Statements a) and b) follow from  the inequalities: $0\leq a \leq m-3$ and  $0\leq b\leq n-3$, combined with the relations 
$$
a=p_{1,j'}-p_{1,j}; \quad b=p_{2,j'}-p_{2,j},
$$

Let us show that Statement c) is true. We have that 
$$
n(m-p_{1,j'}+p_{1,j}-2)+m( n-p_{2,j'}+p_{2,j}-2)= 2nm-na-mb-2n-2m.
$$
We conclude that $2nm-na-mb-2n-2m>nm$ by  Equation \eqref{eq:lema1}.
\end{proof}

\begin{remark}\label{rem:parte inicial derivaciones}
Take a 1-form $\omega$, such that $\nu_C(\omega)=\lambda\notin \Gamma_C$. By Lemma \ref{lem:pr:valores no semi}, we know that $\lambda=j+n+m$ for $j$ an element of the cuspidal set $J$, where we can write $j=np_{1,j}+mp_{2,j}-nm$ with $(p_{1,j},p_{2,j})\in P$. By Proposition \ref{prop:impli:multiplicidad}, we have that $\nu_C(\omega)=n(a+1)+m(b+1)-nm$, where $(a,b)$ is the leading power of a final reduction of $X_\omega(f)$ modulo ${\{f\}}$. Therefore, we conclude that $(a,b)=(p_{1,j},p_{2,j})$. 

In order words, if we have a 1-form $\omega$ whose differential value does not belong to the semigroup, then, the leading power of a final reduction of $X_\omega(f)$ modulo $\{f\}$ is an element of the cuspidal set $P$.
\end{remark}

\subsection{Zariski's invariant case}

In the subsection we proof Theorem \ref{thm:bern:l1}. We see that the theorem is a consequence of the next lemma and a direct application of Theorem \ref{thm:roots}.
\begin{lemma}\label{lem:roots:lambda1}
Denote by $(\lambda_{-1},\lambda_0,\lambda_1,\ldots,\lambda_s)$ the basis of the semimodule of differential values $\Lambda_C$ of $C$. Put $f$ a nice equation of $C$ as in Equation \eqref{eq:pr:nice-equation}. They are equivalent:
\begin{enumerate}
\item $\lambda_1=j_1+n+m$ with $j_1\in J$ is the $1$-element of the basis of $\Lambda_C$.
\item $z_\ell=0$ for $\ell<j_1$ and $z_{j_1}\neq 0$, where the $z_j$ denotes the coefficient of the nice equation $f$. 
\item $Res_f(1,1)((\ell+n+m)/nm)=0$ for $\ell<j_1$ and $Res_f(1,1)((j_1+n+m)/nm)\neq 0$. 
\end{enumerate}
Moreover, we have the following.
\begin{itemize}
\item[$\dagger$)] Assuming that one of the previous conditions is satisfied. For any $\lambda=\lambda_1+na+mb\in \Lambda_C\setminus\Lambda_0$, we have that $Res_f(a+1,b+1)(\lambda/nm)\neq 0$.
\end{itemize}
\end{lemma}

Before proving the Lemma \ref{lem:roots:lambda1}, we remark that the equivalence of the Statements $\mathit 1$ and $\mathit 2$ is well known, see \cite{casas}. However, this equivalence can be proven in an approach similar, but easier, to the one we are going to use in the following section when proving Theorem \ref{thm:bern:4}. Because of that we include our version here. 

\begin{notation}
From now on, given  $r,g\in\mathbb C\{x,y\}$, when we say that $r$ is a reduction of $g$, we mean that $r$ is a reduction of $g$ modulo $\{f\}$. We do similarly with final and partial reductions.
\end{notation}

\begin{proof}[Proof Lemma \ref{lem:roots:lambda1}] We split the proof in three parts: in the first two we prove the equivalence of the three statements, and the third one we prove the Statement $\dagger$).

\bigskip

\emph{Part 1} Statement $\mathit 1$ is equivalent to Statement $\mathit 2$:

\bigskip

Assume that the $1$-element of the basis of the semimodule of differential values of $C$ is $\lambda_1=j_1+na+mb$, with $j_1\in J$. We are going to compute the first elements of a minimal standard basis of the module of differentials of $C$. To do that, we apply Delorme's algorithm.

Since we are in a local system of nice coordinates $(x,y)$, we can put $\omega_{-1}=dx$ and $\omega_0=dy$. Then the axis $u_1$ is 
$$
u_1=\min(\Lambda_0\cap (\lambda_{-1}+\Gamma_C))=n+m=\nu_C(x\omega_0)=\nu_C(y\omega_{-1}),
$$ 
where $\Lambda_i$ denotes the $i^{th}$ element of the decomposition sequence of $\Lambda_C$, with $i=-1,0,\ldots,s$.

We have to consider the 1-form $\theta=x\omega_0+\mu^+y\omega_{-1}$. We need to compute the tuning constant $\mu^+$. Later on, we must find a final reduction of $\theta$ modulo $\{\omega_{-1},\omega_0\}$.

As explained in Section \ref{sec:impli}, by Proposition \ref{prop:impli:multiplicidad}, in order to find $\mu^+$, we need to compute final reductions of $X_{x\omega_0}(f)$ and $X_{y\omega_{-1}}(f)$. Starting by $X_{x\omega_0}(f)$, we find:
\begin{equation}\label{eq:roots: r0}
r_0=X_{x\omega_0}(f)=x\partial_x(f)=xf_x=mx^m+\sum_{\ell\in J}p_{1,\ell}z_\ell x^{p_{1,\ell}}y^{p_{2,\ell}}.
\end{equation}

The leading power $lp(r_0)=(m,0)$ is not divisible by $(0,n)=lp(f)$. Thus, $r_0$ is its own final reduction. For $X_{y\omega_{-1}}(f)$, we have: 
$$
X_{y\omega_{-1}}(f)=-y\partial_y(f)=-yf_y=-ny^n-\sum_{\ell\in J}p_{2,\ell}z_\ell x^{p_{1,\ell}}y^{p_{2,\ell}}.
$$
Here, $lp(X_{y\omega_{-1}}(f))=(0,n)$, so a reduction can be taken as:
\begin{equation}\label{eq:roots: r-1}
r_{-1}=X_{y\omega_{-1}}(f)+nf=nx^m+\sum_{\ell\in J}(n-p_{2,\ell})z_\ell x^{p_{1,\ell}}y^{p_{2,\ell}}.
\end{equation}
Since $lp(r_{-1})=(m,0)$ is non divisible by $(0,n)$, we can put $r_{-1}$ as final reduction  of $X_{y\omega_{-1}}(f)$. By Proposition \ref{prop:pr:codimension}, we find:
$$
i_\mathbf 0(X_{y\omega_{-1}}(f),f)=nm=i_\mathbf 0(X_{x\omega_{0}}(f),f)
$$
From Equations \eqref{eq:roots: r0} and \eqref{eq:roots: r-1}, the leading terms are $lt(r_{-1})=nx^m$ and $lt(r_0)=mx^m$. Therefore, the tuning constant $\mu^+$ is $-m/n$. For convenience, instead of $\theta$, we take the 1-form
$$
\eta=n\theta=nx\omega_{0}-my\omega_{-1}=nxdy-mydx.
$$
We define $r_1$ to be the following partial reduction of $X_{\eta}(f)$:
\begin{equation}\label{eq:roots:r1 lema raices}
\begin{split}
r_1:=nr_0-mr_{-1}&=\sum_{\ell\in J}(np_{1,\ell}+mp_{2,\ell}-nm) z_\ell x^{p_{1,\ell}}y^{p_{2,\ell}}=\\
&=\sum_{\ell\in J}\ell z_\ell x^{p_{1,\ell}}y^{p_{2,\ell}}.
\end{split}
\end{equation}
The leading power of $r_1$ is an element $(p_{1,j},p_{2,j})\in P$. By definition of the cuspidal set $P$, we have that $0\leq p_{2,j}<n-1$. Thus $ (p_{1,j},p_{2,j})$ is not divisible by $(0,n)$. Therefore $r_1$ is a final reduction of $X_\eta(f)$. Moreover, by Propositions \ref{prop:impli:multiplicidad}
\begin{eqnarray*}
i_0(X_\eta(f),f)&=&np_{1,j}+mp_{2,j},\\ \nu_C(\eta)&=&np_{1,j}+mp_{2,j}-(n-1)(m-1)+1=j+n+m.
\end{eqnarray*}

By Lemma \ref{lem:pr:valores no semi}, we have that $j+n+m\notin\Gamma_C=\Lambda_0\cup \{0\}$. We also note that the monomial value of $\eta$ is $\nu_{n,m}(\eta)=n+m=t_1$. Since $\lambda_1$ is the 1-element of the basis, we have that 
$$
\lambda_1=\min(\Lambda_C\setminus \Lambda_0),
$$ 
see Remark \ref{rem:pr:base}. Furthermore, by Theorem \ref{thm:standard:delorme}, it is also satisfied that 
$$
\lambda_1=\sup\{\nu_C(\omega):\nu_{n,m}(\omega)=t_1\}.
$$
We conclude that $\lambda_1=j+n+m$, that is, $j=j_1$. In fact, we have shown that $\nu_C(\eta)=\lambda_1$. This implies that $\eta$ is its own final reduction modulo $\{\omega_{-1},\omega_0\}$, where $r_1$ is a final reduction of $X_\eta(f)$.

By Equation \eqref{eq:roots:r1 lema raices}, stating that the leading power of $r_1$ is $(p_{1,j_1},p_{2,j_1})$ is equivalent to saying that $z_\ell=0$ for $\ell<j_1$ and $z_{j_1}\neq 0$. This shows that Statement $\mathit{1}$ is equivalent to Statement $\mathit{2}$.

\bigskip

\emph{Part 2} Statement $\mathit{2}$ is equivalent to Statement $\mathit{3}$:

\bigskip

Fix $j_1\in J$. Let us compute $Res_f(1,1)((j+n+m)/nm)$, for $j\in J$ with $j\leq j_1$. Recall that $\beta_j=(j+n+m)/nm$, and since $(1,1)\in M$, we can apply Equation  \eqref{eq:residuos}. Notice also that $j=\beta_{j}nm-n-m$, where $\beta_j=(j+n+m)/nm$.

Start by taking $\ell_1=\min(J)$. By Equation \eqref{eq:residuos}, we have to find sequences $(\delta_\ell)_{\ell\in J}$ of non negative integer numbers, such that $\sum_{\ell\in J} \delta_\ell \ell=\ell_1$. Since $\ell_1=\min (J)$, the only possible sequence is the one defined by $\delta_{\ell_1}=1$ and $\delta_\ell=0$ for $\ell\neq \ell_1$. Therefore, by Equation \eqref{eq:residuos}:
$$
Res_f(1,1)((\ell_1+n+m)/nm)=0\Leftrightarrow z_{\ell_1}=0.
$$ 
This proves that $\mathit{2}$ is equivalent to $\mathit{3}$ if $\ell_1=j_1$. 

Now assume that $\ell_1<j_1$, and proceed in an inductive way. Take $k\in J$ such that $k<j_1$. Suppose that $Res_f(1,1)((\ell+n+m)/nm)=0$ for all $\ell\leq k<j_1$ is equivalent to saying that $z_\ell=0$ for all $\ell\leq k$. Denote by $\ell_k=\min\{\ell \in J: \ell >k\}$, let us show that if $z_\ell=0$ for $\ell\leq k$, then
$$
Res_f(1,1)((\ell_k+n+m)/nm)=0\Leftrightarrow z_{\ell_k}=0.
$$
Applying this argument inductively will prove the equivalence between $\mathit{2}$ and $\mathit{3}$. We compute the sequences of non negative integer numbers $(\delta_\ell)_{\ell\in J}$ such that $\sum_{\ell\in J}\delta_\ell\ell=\ell_k$.  There are two kind of possible sequences: first, the one given by $\delta_{\ell_k}=1$ and $\delta_\ell=0$ if $\ell\neq \ell_k$. Second, all the non zero $\delta_\ell$ satisfies that $\ell<\ell_k$. Since $z_\ell=0$ for $\ell<\ell_k$, we have, again by Equation \eqref{eq:residuos}, that
$$
Res_f(1,1)((\ell_k+n+m)/nm)=0\Leftrightarrow z_{\ell_k}=0,
$$
as desired.

\bigskip

\emph{Part 3} Statement $\dagger)$:

\bigskip

Put $\lambda_1=j_1+n+m$ and take $\lambda=\lambda_1+na+mb\in\Lambda_C\setminus\Lambda_0$. By Lemma \ref{lem:pr:valores no semi 2}, we can write $\lambda=q+n+m$ with $q=np_{1,q}+mp_{2,q}-nm$ and $q\in J$. Note that we have the equalities: $p_{1,q}=p_{1,j_1}+a$ and $p_{2,q}=p_{2,j_1}+b$.

As before, we can express:
$$
\frac{q+n+m}{nm}=\frac{np_{1,q}+mp_{2,q}+n+m}{nm}-1=\beta_{q}.
$$
By Lemma \ref{lem:pr:valores no semi 2}, we know that $(a+1,b+1)\in M$. Observe that $j_1=\beta_qnm-n(a+1)-m(b+1).$ Again, since $z_\ell=0$ for $\ell<j_1$, we only have to consider a single non zero sequence $(\delta_\ell)_{\ell\in J}$: $\delta_{j_1}=1$ and $\delta_\ell=0$ for $\ell\neq j_1$. Because, when applying Equation \eqref{eq:residuos} to compute $Res_f(a+1,b+1)(\lambda/nm)$, the other sequences have zero contribution to the residue, thus,
$$
Res_f(a+1,b+1)(\lambda/nm)=0\Leftrightarrow z_{j_1}=0.
$$
However, we have assumed that $z_{j_1}\neq 0$, and we are done.
\end{proof}


\section{Multiplicity up to 4}

As in the previous section, we fix $C$ a cusp with Puiseux pair $(n,m)$, semigroup $\Gamma_C$ and semimodule of differential values $\Lambda_C$. In this section we prove Theorem \ref{thm:bern:4}. Hence, we impose the extra condition that $n\leq 4$. We continue to denote by $\preceq$ the weighted monomial order with respect the weights $(n,m)$. Moreover, we assume that $(x,y)$ is a system of nice coordinates with respect to $C$.  

We recall that the length $s$ of the basis $(\lambda_{-1},\lambda_0,\lambda_1,\ldots,\lambda_s)$ of the semimodule of differential values of $C$ is bounded above by $n-2$. Additionally, we have that $\lambda_{-1}=n$ and $\lambda_0=m$. Therefore, Theorem \ref{thm:bern:4} is trivial if $n=2$. In that case, we have that $\Lambda_C\setminus \Gamma_C=\emptyset$ and there is nothing to proof. By the same argument if $n=3$, we have that $\Lambda_C\setminus \Gamma_C=(\lambda_1+\Gamma_C)\setminus \Gamma_C$. Thus by Theorem \ref{thm:bern:l1}, Theorem \ref{thm:bern:4} is also true when $n=3$. Therefore, we are left to show that it also holds when $n=4$. The rest of the section devoted to show that the theorem holds under the assumption $n=4$. Consider $f\in\mathbb{C}\{x,y\}$ a nice equation of $C$ as in Equation \eqref{eq:pr:nice-equation}:
$$
f=x^m+y^4+\sum_{j\in J} z_j x^{p_{1,j}}y^{p_{2,j}};\qquad z_j\in \mathbb C.
$$
We proceed in a similar way as proving Lemma \ref{lem:roots:lambda1}. We are going to find conditions on the coefficients $z_j$ so that $C$ has $\Lambda_C$ as semimodule of differential values. Later on, we will see that these conditions on the coefficients $z_j$ imply that certain residues, described in Equation \eqref{eq:residuos}, are non zero. Before all of that, we need to determine all the possibilities for $\Lambda_C$.

Since $gcd(4,m)=1$, it follows that $m=4\alpha+\epsilon$ with $\alpha\geq 1$ and $\epsilon\in\{1,3\}$. Before studying all possible semimodules of differential values, we give the following remark about cuspidal sets and nice equations.

\begin{remark}\label{rem: j de 4}
Take $n=4$ and $m=4\alpha+\epsilon$ with $\alpha\geq 2$ and $\epsilon\in\{1,3\}$. For $0\leq \beta \leq \alpha-2$ and $0\leq \beta'\leq 2\alpha-2$ we have that
\begin{eqnarray*}
\epsilon+4\beta &=& 4(3\alpha+\epsilon+\beta)+(4\alpha+\epsilon)1-nm,\\
 2\epsilon+4\beta &=& 4(2\alpha+\epsilon+\beta')+(4\alpha+\epsilon)2-nm.
\end{eqnarray*}
Thus, the cuspidal sets $J$ and $P$ are
$$
J=\{\epsilon+4\beta,\; 2\epsilon+4\beta':\;0\leq \beta\leq \alpha-2,\; 0\leq \beta'\leq 2\alpha-2\}, 
$$
and
$$
P=\{(3\alpha+\epsilon+\beta,1), (2\alpha+\epsilon+\beta',2):\; 0\leq \beta\leq \alpha-2,\; 0\leq \beta'\leq 2\alpha-2\},
$$
with the natural correspondence between $J$ and $P$.

Therefore , we can write the nice equation $f$ of $C$ as:
\begin{equation}\label{eq:4:nice}
f=x^{4\alpha+\epsilon}+y^4+\sum_{\beta=0}^{\alpha-2}z_{\epsilon+4\beta}x^{3\alpha+\epsilon+\beta}y + \sum_{\beta'=0}^{2\alpha-2}z_{2\epsilon+4\beta'}x^{2\alpha+2\epsilon+\beta'}y^2. 
\end{equation}
Observe that, for any $\beta\geq 0$, we have that
$$
(3\alpha+\epsilon+\beta,1)\prec (2\alpha+2\epsilon+\beta,2)	  \prec (3\alpha+\epsilon+\beta+1,1).
$$
\end{remark}

\begin{lemma}\label{lem:4:semimodulos}
Denote by $B$ the basis of the semimodule $\Lambda_C$. Then one of the following statements must be satisfied:
\begin{enumerate}
\item $B=(4,4\alpha+\epsilon)$.
\item $B=(4,4\alpha+\epsilon,\lambda_1)$ with $\lambda_1>u_1=n+m=4(\alpha+1)+\epsilon$ and $\lambda_1\notin \Gamma_C$.
\item $B=(4,4\alpha+\epsilon,\lambda_1,\lambda_2)$ with
\begin{eqnarray*}
\lambda_1&=&n+m+\epsilon+4q=4(\alpha+1)+2\epsilon+4q\quad \text{ with }\quad 0\leq q\leq \alpha-2.\\
\lambda_2&=&8\alpha+3\epsilon+4q'\quad \text{ with }\quad 0\leq q'\leq q.
\end{eqnarray*}
If $\alpha=1$, then Case $\mathit{3}$ is not possible.
\end{enumerate} 
\end{lemma}

Lemma \ref{lem:4:semimodulos} was proven in \cite{hefez3}. There, the authors use their normal form theorem from \cite{hefez2} to show the result. For completeness, we provide a proof using only combinatorial techniques.

\begin{proof}

Put $\lambda_{-1}=n=4$ and $\lambda_0=m=4\alpha+\epsilon$. We recall that, as explained in Section \ref{subsec:cusps}, the only conditions that an increasing sequence $(\lambda_{-1},\lambda_0,\lambda_1,\ldots,\lambda_s)$ must satisfy in order to be the basis of a semimodule of differential values are the following: $\lambda_{-1}=n$, $\lambda_0=m$, $\lambda_i>u_i$ for $i=1,\ldots,s$ and $\lambda_i\notin \lambda_j+\Gamma_C$, for $i\neq j$ and $i,j=-1,0,1,\ldots,s$. It can be verified that the three options above meet these conditions. We are left to show that there are no other possibilities.

We cannot have that the number of the elements of the basis is greater than 4, because $s\leq 2$. Additionally, we see that all possible semimodules of differentials values with $s=0,1$ are covered. For the case $s=1$, we only demand the 1-element of the basis $\lambda_1$ to be greater than $u_1$. Therefore, we must show that if the length of the basis $s$ is 2, the basis must fall into Case $\mathit 3$. 

Hence, assume that $\alpha\geq 2$ and $s=2$; at the end we will show that $\alpha=1$ and $s=2$ are incompatible conditions. In principle there are two possibilities for the 1-element of the basis $\lambda_1$: either $\lambda_1$ is congruent modulo 4 with $2m$, or it is congruent with $3m$.  If $\lambda_1\equiv 3m$ mod 4, let us see that the 2-element of the basis cannot exist. To do that, let us compute the axis $u_2$ and see that for any $\lambda>u_2$, we have that 
$$
\lambda\in \Lambda_1=(n+\Gamma_C)\cup (m+\Gamma_C)\cup (\lambda_1+\Gamma_C).
$$  
By Remark \ref{rem:pr:escritura ui}, we know that $u_2$ must take one of the following four values:
\begin{eqnarray*}
\lambda_1+mp&=\hphantom{ii} n+na; \qquad \lambda_1+mp'&=m+na';\\
\lambda_1+nq&=m+mb;\phantom{h} \qquad \lambda_1+nq'&=n+mb';
\end{eqnarray*}
where $p,p'q,q'$ and $a,a',b,b'$ are the minimal non negative integer numbers satisfying the previous equations. In particular, the axis $u_2$ is the smallest of the previous four values. Note that $p'=p+1$ and $q'=q+1$, thus the cases from the right column are discarded.  

If $u_2=\lambda_1+mp=n(a+1)$, since $n=4$ and $\lambda_1\equiv 3m$ mod 4, then, $p=1$. Therefore, for any integer $\lambda>u_2$, we find that $\lambda>u_2=\lambda_1+m>2m$. Assume that we write $\lambda\equiv \delta m$ mod 4 with $\delta\in \{0,1,2,3\}$. If $\delta=0,1,2$, we could write $\lambda=\delta m+4c$ for $c\geq 0$ and $\lambda_2\in \Gamma_C$. If $\delta=3$, then $\lambda=\lambda_1+4c'$ with $c'\geq 0$ and $\lambda\in \lambda_1+\Gamma_C$. In both cases, we have that $\lambda\in \Lambda_1$

If $u_2=\lambda_1+nq=m(b+1)$, we have that $b\geq 2$, thus $\lambda>u_2\geq 3m$. As before, we could write $\lambda=\delta m+4c$ with $c\geq 0$ and $\delta\in \{0,1,2,3\}$. This means that $\lambda\in\Gamma_C$. In other words, if $\lambda_1\equiv 3m$ mod 4, we obtain that $\lambda>u_2$ implies $\lambda\in\Lambda_1$. 

In conclusion the assumption $\lambda_1\equiv 3m$ mod 4 implies that the 2-element of the basis cannot exist.

Now assume that $\lambda_1\equiv 2m$ mod 4. 
First, we note that $\lambda_1<2m$, otherwise, we could write $\lambda_1=2m+4\ell$ for $\ell\geq 0$ implying that $\lambda_1\in \Gamma_C$. This would contradict the fact that $\lambda_1$ is an element of the basis different from $n$ and $m$. Moreover, we have the extra condition $\lambda_1>u_1=4(\alpha+1)+\epsilon$. In other words 
$$
4(\alpha+1)+\epsilon<\lambda_1<2m=4(2\alpha)+2\epsilon.
$$
Taking into account that we are assuming that $\lambda_1\equiv 2m$ mod 4, the last two inequalities are equivalent to 
$$
\lambda_1=4(\alpha+1)+2\epsilon+4q,\quad \text{ with }\quad 0\leq q\leq \alpha-2.
$$
Hence $\lambda_1$ must be as stated in Case \textit{3}. We are left to determine all the possibilities for the 2-element of the basis. We can check that $u_2=\lambda_1+4(\alpha-q-1)=2m=8\alpha+2\epsilon$. Assume that $\lambda>u_2$ and that $\lambda\notin \Lambda_1$, as the 2-element of the basis must satisfy.

The conditions $\lambda>u_2>\lambda_1$ and  $\lambda\notin \Lambda_1$ imply that $\lambda\equiv 3m$ mod 4. Otherwise, we could write $\lambda=\lambda_k+4c$ for $k=-1,0,1$ and $c\geq 0$, where, as always, $\lambda_{-1}=n$ and $\lambda_0=m$. Arguing as before, we find that 
$$
u_2<\lambda<\lambda_1+m=8\alpha+4+4q+3\epsilon <3m.
$$
From the previous inequalities and $\lambda\equiv 3m$ mod 4 we get 
$$
\lambda=8\alpha+3\epsilon+4q'\quad \text{ with }\quad 0\leq q'\leq q.
$$
This shows that the basis is as in Case $\mathit 3$.

\bigskip

Finally, we notice that if $\alpha=1$, it is not possible to have a cusps whose semimodules of differential values whose basis has length 2. Indeed, assume that $\epsilon=1$, then we have that $\Gamma_C=\langle 4,5\rangle$. The conductor of the semigroup is $c_\Gamma=(4-1)(5-1)=12$. Then the axis that we compute is $u_1=9$. This implies that if the 1-element of the basis $\lambda_1$ exists, then it must be $\lambda_1=11>9$ (note that $10=2\cdot5\in \Gamma_C$). Thus, we compute the new axis and we have that $u_2=15>c_\Gamma$. Therefore given $\lambda>u_2$, we see that $\lambda\in \Gamma_C$. We conclude that the semimodule of differential values of a cusp with Puiseux pair $(4,5)$ cannot have a basis of length 2. If $\epsilon=3$ we proceed in a similar way.
\end{proof}

\begin{remark}\label{rem:4:elementos semimodulo}
The proof of Lemma \ref{lem:4:semimodulos} shows that in case of having a basis as in Case $\mathit 3$, then for any $\lambda\notin \Lambda_1$ with $\lambda>u_2$, we have that 
$$
\lambda=8\alpha+3\epsilon+4\beta\quad \text{ with }\quad 0\leq \beta\leq q.
$$
Now, consider the 2-element of the basis $\lambda_2=8\alpha+3\epsilon+4q'$ with $0\leq \beta\leq q'$. Given $\lambda\in\Lambda_2\setminus \Lambda_1=\Lambda_C\setminus \Lambda_1$, then $\lambda=8\alpha+3\epsilon+4\beta$ with $q'\leq \beta\leq q$. This is equivalent to saying that $\lambda=\lambda_2+na$ with $0\leq a\leq q-q'$.
\end{remark}

By Lemma \ref{lem:4:semimodulos}, we only have to prove the Theorem \ref{thm:bern:4} in three cases. In the first one, where the length of the basis is $s=0$, there is nothing to prove because $\Lambda_C\setminus \Gamma_C=\emptyset$. In the second one, with $s=1$, since $\Lambda_C\setminus \Gamma_C=(\lambda_1+\Gamma_C)\setminus \Gamma_C$, the theorem holds by Theorem \ref{thm:bern:l1}. Hence we are left to show the third case with $s=2$. In that case, again by Theorem \ref{thm:bern:l1}, we only need to show that any element in $\lambda\in \Lambda_C\setminus \Lambda_1$, then $-\lambda/nm$ becomes a root of the Bernstein-Sato polynomial. Moreover the elements $\Lambda_C\setminus \Lambda_1$ are described in the previous remark.

The proof of Theorem \ref{thm:bern:4} for this last case is straight forward by the next lemma and Theorem \ref{thm:roots}.

\begin{lemma}\label{lem:4:equiv1}
With the previous notations, assume that $(n=4,m=4\alpha+\epsilon,\lambda_1=4(\alpha+1)+2\epsilon+4q,\lambda_2)$ is the basis of the semimodule of differential values of $C$,  where $\alpha\geq 2$, $\epsilon\in	\{1,3\}$ and $0\leq q\leq \alpha-2$. They are equivalent:
\begin{enumerate}
\item $\lambda_2=8\alpha+3\epsilon+4q'$ with $0\leq q'\leq q$ is the 2-element of the basis.
\item $Res_f(\alpha-q,1)((8\alpha+3\epsilon+4\gamma)/nm)=0$ for all non negative integers $\gamma<q'$ and $Res_f(\alpha-q,1)((8\alpha+3\epsilon+4q')/nm)\neq 0$.
\end{enumerate}
Moreover, it also holds the following statement:
\begin{itemize}
\item[$\dagger)$] Assuming one of the two above conditions are satisfied, for any $\lambda_2+na+mb\notin \Lambda_1$, we have that  $Res_f(\alpha-q+a,b+1)((\lambda_2+na+mb)/nm)\neq 0$.
\end{itemize}
\end{lemma}

As in Lemma \ref{lem:roots:lambda1}, the first two statements are equivalent to some algebraic conditions on the coefficients $z_j$ of the equation $f=0$. 

\begin{lemma}\label{lem:4:equiv2}
With the previous notations, assume that $(n=4,m=4\alpha+\epsilon,\lambda_1=4(\alpha+1)+2\epsilon+4q,\lambda_2)$ is the basis of the semimodule of differential values of $C$,  where $\alpha\geq 2$, $\epsilon\in	\{1,3\}$ and $0\leq q\leq \alpha-2$. They are equivalent:
\begin{enumerate}
\item $\lambda_2=8\alpha+3\epsilon+4q'$ with $0\leq q'\leq q$ is the 2-element of the basis.
\item $z_{2\epsilon+4\beta'}=0$ for $q\leq \beta' < q+q'$ and $z_{2\epsilon+4(q+q')}\neq 0$, if $q'<q$. Or 
$$
2(4\alpha+\epsilon) z_{2\epsilon+8q}-(3\alpha+\epsilon+q)z^2_{\epsilon+4q}\neq 0,
$$
if $q'=q$. 
\end{enumerate}
\end{lemma}

Before giving the proof, we recall that the term ``reduction'' refers specifically to ``reduction modulo $\{f\}$'' when referring to functions.

\begin{proof}[Proof Lemma \ref{lem:4:equiv2}]

Assume that the 2-element of the basis of $\Lambda_C$ is $\lambda_2=8\alpha+3\epsilon+4q'$ with $0\leq q'\leq q$. Essentially, the proof follows the same reasoning as in Lemma \ref{lem:roots:lambda1} part 1. We will apply the Delorme's algorithm to compute a minimal standard basis of $C$. This way we will obtain a 1-form $\omega_2$ whose differential value is $\lambda_2$. During the process, we will derive the desired algebraic conditions from Statement \textit{2}.


In order to simplify the computations along all the proof, we use the following notation: consider $g\in\mathbb{C}\{x,y\}$ and assume that $g=\sum a_{jk}x^jy^k$. When we write
$$
g=g_1+h.t.,
$$
we mean that $g_1=\sum b_{jk}x^jy^k$ where $b_{jk}=0$ for $(j,k)\succ(3\alpha+\epsilon-1+q,2)$ and $b_{jk}=a_{jk}$ otherwise. Additionally, the higher terms $(h.t.)$ are those monomials $a_{jk}x^jy^k$ of $g$ with $(j,k)\succ (3\alpha+\epsilon-1+q,2)$. 

We introduce the previous notation for the following reason: consider a 1-form $\eta$ with monomial value $\nu_{n,m}(\eta)=t_2$. By Theorem \ref{thm:standard:delorme}, we know that $\nu_C(\eta)\leq \lambda_2$. Thus, by Propositions  \ref{prop:impli:base-ideal} and \ref{prop:impli:multiplicidad}, we have that
$$
i_\mathbf 0(X_\eta(f),f)\leq \lambda_2+nm-n-m=4(3\alpha+q'+\epsilon-1)+(4\alpha+\epsilon)2.
$$
Therefore, if $r$ is a partial or final reduction of $X_\eta(f)$, then $r$ satisfies the following property:
\begin{itemize}
\item [$\star$)] The leading power of $r$ is at most $(3\alpha+q'+\epsilon-1,2)$. 
\end{itemize}
Since $q'\leq q$, we are not concerned with the behaviour of the monomials with leading power greater than $(3\alpha+q+\epsilon-1,2)$.

%
%

That been said, we start by computing a 1-form that could later be identified as $\omega_2$. 
Since $\lambda_1=4(\alpha+1)+2\epsilon+4q=m+n+\epsilon+4q$ with $0\leq q\leq \alpha-2$, by Lemma \ref{lem:roots:lambda1}, we have that $z_\ell=0$ for $\ell<\epsilon +4q$ and $z_{\epsilon+4q}\neq 0$. Hence, we can change the index of both summations in the Equation \eqref{eq:4:nice} as
\begin{equation}\label{eq:nice lema reducida}
f=x^{4\alpha+\epsilon}+y^4+\sum_{\beta=q}^{\alpha-2}z_{\epsilon+4\beta}x^{3\alpha+\epsilon+\beta}y + \sum_{\beta'=q}^{2q}z_{2\epsilon+4\beta'}x^{2\alpha+\epsilon+\beta'}y^2.
\end{equation}
We now apply the Delorme's algorithm. Since we are in a system of nice coordinates $(x,y)$ to $C$, if we set $\omega_{-1}=dx$ and $\omega_0=dy$, then $\nu_C(\omega_{-1})=n$ and $\nu_C(\omega_0)=m$. Similarly to the proof of Lemma \ref{lem:roots:lambda1}, we find that the 1-form $\omega_1=nxdy-mydx$ satisfies that the function $r_1$ given by
\begin{equation}\label{eq:lema r1}
\begin{split}
r_1&=X_\omega(f)-nmf=nxf_x+myf_y-nmf=\sum_{j\in J} jz_jx^{p_{1,j}}y^{p_{2,j}}=\\
&= \sum_{\beta=q}^{\alpha-2}(\epsilon+4\beta)z_{\epsilon+4\beta}x^{3\alpha+\epsilon+\beta}y + \sum_{\beta'=q}^{2q}(2\epsilon+4\beta')z_{2\epsilon+4\beta'}x^{2\alpha+\epsilon+\beta'}y^2,
\end{split}
\end{equation}
it is a final reduction of $X_{\omega_1}(f)$. Moreover, the leading power of $r_1$ is $lp(r_1)=(3\alpha+\epsilon+q,1)$. By Proposition \ref{prop:impli:multiplicidad}, the differential value of $\omega_1$ is $\nu_C(\omega_1)=\lambda_1$.

We now compute our candidate for $\omega_2$. As in the proof of Lemma \ref{lem:4:semimodulos}, the axis $u_2$ is $u_2=\lambda_1+4(\alpha-q-1)=2m$. In particular, we have: 
$$
u_2=\nu_C(x^{\alpha-q-1}\omega_1)=\nu_C(ydy).
$$ 
This implies that  $t_2=t_1+4(\alpha-q-1)=4(\alpha-q)+m<2m$. Therefore, by Delorme's algorithm we need to compute the tuning constant $\mu^+$ from the 1-form $\theta=x^{\alpha-q-1}\omega_1+\mu^+ydy$. Later, we have to compute a final reduction of $\theta$ modulo $\{\omega_{-1},\omega_0,\omega_1\}$.

\bigskip

\emph{Computation of $\mu^+$}: We must compute final reductions of $X_{x^{\alpha-q-1}\omega_1}(f)$ and $X_{ydy}(f)$. For the first one, by Equation \eqref{eq:lema r1}, we have that
\begin{equation}\label{eq:lema r2}
\begin{split}
r_2&=x^{\alpha-q-1}r_1=x^{\alpha-q-1}(X_{\omega_1}(f)-nmf)\\
&=\sum_{\beta=q}^{\alpha-2}(\epsilon+4\beta)z_{\epsilon+4\beta}x^{4\alpha+\epsilon+\beta-q-1}y + \\
&+\sum_{\beta'=q}^{2q}(2\epsilon+4\beta')z_{2\epsilon+4\beta'}x^{3\alpha+\epsilon+\beta'-q-1}y^2,
\end{split}
\end{equation}
which is indeed a final reduction of $X_{x^{\alpha-q-1}\omega_1}(f)$, as its leading power $(4\alpha+\epsilon-1,1)$ is not divisible by $(0,4)$. Similarly, for $X_{ydy}(f)$, we have that
\begin{equation}\label{eq:lema r0}
\begin{split}
r_0&= X_{ydy}(f)= \\
&= (4\alpha+\epsilon)x^{4\alpha+\epsilon-1}y +(3\alpha+\epsilon+q)z_{\epsilon+4q}x^{3\alpha+\epsilon+q-1}y^2 + h.t.
\end{split}
\end{equation}
Again, the leading power $lp(r_0)=(4\alpha+\epsilon-1,1)$ is not divisible by $(0,4)$. Therefore, by Equations \eqref{eq:lema r2} and \eqref{eq:lema r0}, the tuning constant is $\mu^+=-(\epsilon+4q)z_{\epsilon+4q}/(4\alpha+\epsilon)$. Equivalently, we can write: 
\begin{equation}\label{eq:lema w}
\omega= (4\alpha+\epsilon)\theta=(4\alpha+\epsilon)x^{\alpha-q-1}\omega_1-(\epsilon+4q)z_{\epsilon+4q}ydy.
\end{equation}
This way we have that $\nu_C(\omega)>\nu_C(x^{\alpha-q-1}\omega_1)=\nu_C(ydy)$. Notice that $\nu_{n,m}(\omega)=t_2=4(\alpha-q)+m$, because $\nu_{n,m}(x^{\alpha-q-1}\omega_1)=t_2=4(\alpha-q)+m<2m=\nu_{n,m}(ydy)$. 

\bigskip

\emph{Computation of a final reduction of $\omega$}: Note that it is the same than computing a final reduction of $\theta$. For this purpose, we should recursively construct 1-forms $\omega^i$ of the form
$$
\omega^{i}=\omega^{i-1}+\mu^+ h_i\omega_k,
$$ 
where $k\in\{-1,0,1\}$, $i\geq 0$, $h_i$ is a monomial, and $\omega^0=\omega$. This process continues until we obtain a 1-form whose differential value is $\lambda_2$.

To shorten the proof we do it in a single step, we take all the reductions modulo $\{\omega_{-1},\omega_0,\omega_1\}$ simultaneously. To do so, we compute a partial reduction $X$ of $X_\omega(f)$. Using Equations \eqref{eq:lema r2} and \eqref{eq:lema r0}, we define $X$ as
\begin{equation}\label{eq:lema X}
\begin{split}
X& =(4\alpha+\epsilon)r_2-((\epsilon+4q)z_{\epsilon+4q})r_0=\\
&=\sum_{\beta=q+1}^{2q}(4\alpha+\epsilon)(\epsilon+4\beta)z_{\epsilon+4\beta}x^{4\alpha+\epsilon+\beta-q-1}y +\\ 
&+ \sum_{\beta'=q}^{2q}(4\alpha+\epsilon)(2\epsilon+4\beta')z_{2\epsilon+4\beta'}x^{3\alpha+\epsilon+\beta'-q-1}y^2 - \\
&- (\epsilon+4q)(3\alpha+\epsilon+q)z_{\epsilon+4q}^2x^{3\alpha+\epsilon+q-1}y^2+h.t.,
\end{split}
\end{equation}
since $r_0$ and $r_2$ are final reductions of $X_{ydy}(f)$ and $X_{x^{\alpha-q-1}\omega_1}(f)$ respectively. Comparing the Equation \eqref{eq:lema w} with Equation \eqref{eq:lema X}, we verify that $X$ is a partial reduction of $X_\omega(f)$. 

Since $\nu_{n,m}(\omega)=t_2$, then by the $\star)$ property, we have that $lp(X)\preceq (3\alpha+q'+\epsilon-1,2)$. In particular, this implies that $X\neq 0$. Moreover, by Equation \eqref{eq:lema X}, we observe that $lp(X)\succeq (3\alpha+\epsilon-1,2)$. 

%

The function $X$ encodes the necessary information to compute the desired final reduction of $\omega$. Write 
$$
X=\sum_{j,k\geq 0} a_{jk}x^jy^k,
$$
and assume that there exists a  minimum index $q\leq \ell\leq q+q'$, satisfying that $a_{3\alpha+\epsilon+\ell-q-1,2}\neq 0$. Note that by Equation \eqref{eq:lema X}, having $a_{3\alpha+\epsilon+\beta'-q-1,2}\neq 0$ is the same as having $z_{2\epsilon+4\beta'}\neq 0$, where $q\leq \beta'<2q$. In other words, we are assuming that $z_{2\epsilon+4\beta'}=0$ for $q\leq \beta'<\ell$ and that $z_{2\epsilon+4\ell}\neq 0$.

We take the 1-form $\omega'$ which is going to be a final reduction of $\omega$ modulo $\{\omega_{-1},\omega_0,\omega_1\}$. It is given by:
\begin{equation}\label{eq:lema w'}
\omega'=\omega-\sum_{\beta=q+1}^\ell (\epsilon+4\beta)z_{\epsilon+4\beta} x^{\beta-q}ydy.
\end{equation}
Notice that for any $\beta >q$, we have that 
$$
\nu_{n,m}(x^{\beta-q}ydy)>\nu_{n,m}(ydy)=2m>t_2=4(\alpha-q)+m=8\alpha+\epsilon-4q.
$$
Thus, $\nu_{n,m}(\omega')=t_2$. Next, observe that:
$$
X_{x^{\beta-q}ydy}(f)=(4\alpha+\epsilon)x^{4\alpha+\epsilon+\beta-q-1}y+h.t.,
$$
where we see that $X_{x^{\beta-q}ydy}(f)$ is non reducible modulo $\{f\}$. Since $X$ is a partial reduction of $X_\omega(f)$, we can define a partial reduction $X'$ of $X_{\omega'}(f)$, given by the expression:
\begin{equation}\label{eq:lema X'}
\begin{split}
X'&=X-\sum_{\beta=q+1}^\ell (\epsilon+4\beta)z_{\epsilon+4\beta} X_{x^{\beta-q}ydy}(f)=\\
&=\sum_{\beta=\ell+1}^{2q}(4\alpha+\epsilon)(\epsilon+4\beta)z_{\epsilon+4\beta}x^{4\alpha+\epsilon+\beta-q-1}y +\\ 
&+ \sum_{\beta'=\ell}^{2q}(4\alpha+\epsilon)(2\epsilon+4\beta')z_{2\epsilon+4\beta'}x^{3\alpha+\epsilon+\beta'-q-1}y^2 - \\
&- (\epsilon+4q)(3\alpha+\epsilon+q)z_{\epsilon+4q}^2x^{3\alpha+\epsilon+q-1}y^2+h.t.
\end{split}
\end{equation}
Since $z_{2\epsilon+4\ell}\neq 0$, we have that the leading power of $X'$ is $lp(X')=(3\alpha+\epsilon+\ell-q-1,2)$, which is not divisible by $lp(f)=(0,4)$. Hence, $X'$ is a final reduction of $X_{\omega'}(f)$. By Proposition \ref{prop:impli:multiplicidad}, we have that 
$$
\nu_C(\omega')=8\alpha+3\epsilon+4(\ell-q)=\lambda.
$$ 
We note the following: first, $\lambda\leq \lambda_2$, with the equality achieved if and only if $\ell=q+q'$. Second, $\lambda\notin \Lambda_1$. Thus, for $\lambda_2$ to be the minimum element in $\Lambda_C\setminus \Lambda_1$, we need $\ell=q+q'$. 

Finally, if $\ell$ does not exist, meaning $a_{3\alpha+\epsilon+\beta'-1,2}= 0$ for $q\leq \beta'\leq q+q'$, we can construct $X'$ as before by setting $\ell=q+q'$ in the expression of $\omega'$ from Equation \eqref{eq:lema w'}. However, this time, by Equation \eqref{eq:lema X'}, $X'$ may not be a final reduction of $X_{\omega'}(f)$. Nonetheless, we see  that its leading power is greater than $(3\alpha+\epsilon+q'-1,2)$. But, as we saw before, that is not possible, since $\nu_{n,m}(\omega')=t_2$ and the $\star)$ property.

Therefore, $\ell$ exists and it takes the value $\ell=q+q'$. In this situation, we put $\omega_2=\omega'$ and $X'$ is a final reduction of $X_{\omega_2}(f)$. 

\bigskip

\emph{Conclusion}: As mentioned before, by Equation \eqref{eq:lema X'}, having $a_{3\alpha+\epsilon+\beta'-q-1,2}=0$ for $q\leq \beta'<q+q'$ is equivalent to having $z_{2\epsilon+4\beta'}=0$. Additionally, the condition $a_{3\alpha+\epsilon+q'-1,2}\neq 0$ is the same as $z_{2\epsilon+4(q+q')}\neq 0$ if $q'<q$, and to 
\begin{eqnarray*}
(4\alpha+\epsilon)(2\epsilon+8q)z_{2\epsilon+8q}-(\epsilon+4q)(3\alpha+\epsilon+q)z_{\epsilon+4q}^2&\neq &0 \Leftrightarrow\\
2(4\alpha+\epsilon)z_{2\epsilon+8q}-(3\alpha+\epsilon+q)z_{\epsilon+4q}^2&\neq &0,
\end{eqnarray*}
if $q'=q$, as desired. Finally, all the previous computations also shows that if we assume Statement \textit{2}, then $\nu_C(\omega_2)=8\alpha+3\epsilon+4q'$. 

\end{proof}

Now, we give the proof of Lemma \ref{lem:4:equiv1}.

\begin{proof}[Proof Lemma \ref{lem:4:equiv1}]

We split the proof into two parts: in the first one, we show the equivalence between the two results. In the second part, we prove the final statement of the lemma.

\bigskip

\emph{Part 1} Statement \textit{1} is equivalent to Statement \textit{2}:

\bigskip

Assume that $\lambda_2=8\alpha+3\epsilon+4q'$ with $0\leq q'\leq q$. We need to compute the residues $Res_f(\alpha-q,1)((8\alpha+3\epsilon+4\gamma)/nm)$ for $0\leq \gamma\leq q'$. First, we claim that that $(\alpha-\delta,1)$ belongs to the cuspidal set $M$ for $0\leq \delta\leq \alpha-2$, and in particular, $(\alpha-q,1)\in M$, showing that we are in position to apply Equation \eqref{eq:residuos}. To verify this, note that
$$
(\alpha-\delta,1)=(m-1,n-1)-(p_1,p_2),
$$ 
where $(p_1,p_2)=(3\alpha+\epsilon+\delta-1,2)$ is an element of the cuspidal set $P$, as shown in Remark \ref{rem: j de 4}. 

Now, consider $0\leq \gamma< q'\leq q$. We observe that: 
$$
8\alpha+3\epsilon+4\gamma-4(\alpha-q)-(4\alpha+\epsilon)1=2\epsilon+4(q+\gamma)\in J.
$$
Hence, in virtue of Equation \eqref{eq:residuos}, we need to find the sequences of non negative integer numbers $(\delta_\ell)_{\ell\in J}$, such that: 
$$
\sum_{\ell\in J}\ell\delta_\ell=2\epsilon+4(q+\gamma).
$$
Recall from Remark \ref{rem: j de 4} that the elements in $J$ are of the form $\epsilon+4\beta$ and $2\epsilon+4\beta'$, where $0\leq \beta\leq \alpha-2$ and $0\leq \beta'\leq 2\alpha-2$.  As shown in the proof of Lemma \ref{lem:4:equiv2}, we can take a nice equation $f$ of $C$, given by: 
$$
f=x^{4\alpha+\epsilon}+y^4+\sum_{\beta=q}^{\alpha-2}z_{\epsilon+4\beta}x^{3\alpha+\epsilon+\beta}y + \sum_{\beta'=q}^{2\alpha-2}z_{2\epsilon+4\beta'}x^{2\alpha+\epsilon+\beta'}y^2.
$$
Hence, there is only one sequence of $(\delta_\ell)_{\ell\in J}$ that is relevant in the computation of the residue: $\delta_{2\epsilon+4(q+\gamma)}=1$ and $\delta_\ell=0$ for $\ell\neq 2\epsilon+4(q+\gamma)$. 
Since any other sequence $(\delta_\ell')_{\ell\in J}$  have 0 contribution to the computation of the residue, by Equation \eqref{eq:residuos}, we conclude that:
\begin{equation}\label{eq:lema: residuos 1}
Res_f(\alpha-q,1)((8\alpha+3\epsilon+4\gamma)/nm)\neq 0 \Leftrightarrow z_{2\epsilon+4(\gamma+q)}\neq 0.
\end{equation}
Now, if $\gamma=q'$, there are two cases: $q'<q$ or $q'=q$. Assume first that $q'<q$, then we have, as before, a single relevant sequence: $\delta_{2\epsilon+4(q+q')}=1$ and $\delta_\ell=0$ if $\ell\neq 2\epsilon+4(q+q')$. This gives: 
\begin{equation}\label{eq:lema: residuos 2}
Res_f(\alpha-q,1)((8\alpha+3\epsilon+4q')/nm)\neq 0 \Leftrightarrow z_{2\epsilon+4(q+q')}\neq 0.
\end{equation}
If $q'=q$, then there are two relevant sequences: first, $\delta_{2\epsilon+8q}=1$ and $\delta_\ell=0$ if $\ell\neq 2\epsilon+8q$. Second, $\delta_{\epsilon+4q}=2$ and $\delta_\ell=0$ if $\ell\neq \epsilon+4q$. Hence, by Equation \eqref{eq:residuos} we obtain the following:
\begin{eqnarray*}
Res_f(\alpha-q,1)\left(\tfrac{8\alpha+3\epsilon+4q}{nm}\right)&=&\frac{\Gamma\left(\tfrac{8\alpha+3\epsilon+4q}{nm}\right)^{-1}}{nm}\left[\tfrac{z_{\epsilon+4q}^2}{2}\Gamma\left( \tfrac{2(3\alpha+\epsilon+q)+(\alpha-q)}{m}\right)\Gamma\left(\tfrac{3}{4}\right)-\right.\\
&-&\left. z_{2\epsilon+8q}\Gamma\left( \tfrac{(2\alpha+\epsilon+2q)+(\alpha-q)}{m}\right)\Gamma\left(\tfrac{3}{4}\right) \right].
\end{eqnarray*}
Using the fact that the Euler's Gamma function satisfies that $\Gamma(\rho+1)=\rho\Gamma(\rho)$, we can extract a common factor in the previous equation. This leads to:
\begin{equation}\label{eq:lema: residuos 3}
Res_f(\alpha-q,1)\left(\tfrac{8\alpha+3\epsilon+4q}{nm}\right)\neq 0 \Leftrightarrow 2(4\alpha+\epsilon) z_{2\epsilon+q}-(3\alpha+\epsilon+q)z^2_{\epsilon+4q}\neq 0.
\end{equation}

By Equations \eqref{eq:lema: residuos 1}-\eqref{eq:lema: residuos 3}, we have to show that Statement \textit{1} is equivalent to 
\begin{itemize}
\item $z_{2\epsilon+4(\gamma+q)}=0$ for $0\leq \gamma<q'$.
\item if $q'<q$, then $z_{2\epsilon+4(q'+q)}\neq 0$.
\item if $q'=q$, then $z_{2\epsilon+q}-(3\alpha+\epsilon+q)z^2_{\epsilon+4q}\neq 0$
\end{itemize}

Thus, we see that Statement $\mathit 1$ and $\mathit 2$ are equivalent in virtue of Lemma \ref{lem:4:equiv2}.

\bigskip

\emph{Part 2} Statement $\dagger)$:

\bigskip

We need to show that for a given $\lambda=\lambda_2+na+mb\notin\Lambda_1$ with $a,b\geq 0$, then  
$$
Res_f(\alpha-q+a,b+1)(\lambda/nm)\neq 0.
$$
By Remark \ref{rem:4:elementos semimodulo}, since $\lambda\notin \Lambda_1$, we have that $b=0$ and $0\leq a\leq q-q'$. If $a=0$, we have already shown that $Res_f(\alpha-q,1)(\lambda_2/nm)\neq 0$.  Thus, assume that $1\leq a\leq q-q'$, notice that we are assuming that $q'<q$. As shown at the beginning of the part 1 of the proof, we have $(\alpha-q+a,1)\in M$. 

We can write $\lambda$ explicitly as $\lambda=8\alpha+3\epsilon+4q'+4a$. Substracting $n(\alpha-q+a)+m$, we find:  
$$
\lambda-n(\alpha-q+a)-m=2\epsilon+4(q+q').
$$ 
In other words, we have to compute sequences $(\delta_\ell)_{\ell\in J}$ such that $\sum \ell \delta_\ell=2\epsilon+4(q+q')$. Since $q'<q$, there is only one relevant sequence to the computation of the desired residue: $\delta_{2\epsilon+4(q+q')}=1$ and the rest $\delta_\ell=0$ for $\ell\neq 2\epsilon+4(q+q')$. Any other sequence gives zero contribution to the residue. Thus, 
$$
Res_f(\alpha-q+a,1)(\lambda/nm)\neq 0\Leftrightarrow z_{2\epsilon+4(q+q')}\neq 0
$$
and we have shown that $z_{2\epsilon+4(q+q')}\neq 0$.

\end{proof}

\appendix

\section{Minimal Standard basis of the extended Jacobian ideal}\label{app:jac}

Assume that $C$ is a cusp with Puiseux pair $(n,m)$, let $f=0$ be an implicit equation of $C$ and $(x,y)$ is an adapted system of coordinates with respect to $C$. Additionally, consider the weighted order $\preceq$ with respect $(n,m)$. The goal of this appendix is to show how a minimal standard basis of the module of differentials is related with a minimal standard basis of the extended jacobian ideal of $C$, when considering monomial order $\preceq$. Recall that the extended jacobian ideal of $f$ is $\mathcal J(f)=(f,f_x,f_y)$, where $f_x,f_y$ are the partial derivatives of $f$ with respect $x$ and $y$ respectively.

In \cite{brianson}, it is provided a minimal standard basis of the extended jacobian ideal, when $f$ is generic. Besides, according to \cite{pol}, if we have a minimal standard basis of $\mathcal J(f)$, then we can obtain the semimodule of differential values of $C$.

\begin{theorem}\label{thm:jac:base-estandar}
Assume that $C$ is a cusp with a Puiseux pair $(n,m)$ and let $f\in\mathbb C\{x,y\}$ be an implicit equation. Suppose also, that the local system of coordinates $(x,y)$ is adapted with respect $C$. Denote by $(\lambda_{-1},\lambda_0,\ldots,\lambda_s)$ the basis of the semimodule of differential values $\Lambda_C$ of $C$. Take $(\omega_{-1},\omega_0,\ldots,\omega_s)$ a minimal standard basis of the module of differentials of $C$, and for $i=-1,0,\ldots,s$ put $h_i\in\mathbb C\{x,y\}$ a final reduction of $X_{\omega_i}(f)$ modulo $\{f\}$, then 
$$
B=\{h_{-1},h_0,\ldots,h_s\}
$$
is a minimal standard basis of $\mathcal J(f)$ with respect the weighted order with respect $(n,m)$.
\end{theorem}
\begin{proof}
Put
$$
f=\mu x^m+y^n+\sum_{\substack{\alpha,\beta\geq 0\\ n\alpha+m\beta >nm}}z_{\alpha\beta}x^\alpha y^\beta.
$$
It is enough to show the following three statements: 
\begin{enumerate}
\item $B$ is a generator system of the ideal $\mathcal J(f)$.
\item Given $g\in \mathcal J(f)$, there exists at least one element $b\in B$, such that $lp(b)$ divides $lp(g)$.
\item Given $h_i,h_j\in B$ with $h_i\neq h_j$, then the leading powers satisfy that $lp(h_i)\nmid lp(h_j)$.
\end{enumerate}
\emph{Statement} 1: We need to show that the ideal generated by $B$ coincides with $\mathcal J(f)$. Given $\omega\in\Omega_{\mathbb C^2,\mathbf 0}^1$, by definition we have that any final reduction of $X_{\omega}(f)$ is an element of $\mathcal J(f)$. Therefore, we only need to show that $f_x,f_y$ and $f$ belong to the ideal $(B)$.

By Remark \ref{rem:pr:xy}, the pair $\{\omega_{-1},\omega_0\}$ is a basis of the $\mathbb C\{x,y\}$-module $\Omega_{\mathbb C^2,\mathbf 0}^1 $. Hence, we can write the 1-form $dx$ as 
\begin{equation}\label{eq:ap:dx}
dx=A\omega_{-1}+B\omega_0\text{ with }A,B\in \mathbb{C}\{x,y\}.
\end{equation}
Besides, we also have that $\nu_C(\omega_{-1})=n$ and $\nu_C(\omega_0)=m$, thus by Equation \eqref{eq:impli:valor-diferencial}, we see that $i_{\mathbf 0}(X_{\omega_{-1}}(f),f)=nm-m$ and $i_{\mathbf 0}(X_{\omega_{0}}(f),f)=nm-n$. Therefore, by Proposition \ref{prop:impli:multiplicidad}, the leading powers of $h_{-1}$ and $h_0$ are $lp(h_{-1})=(0,n-1)$ and $lp(h_0)=(m-1,0)$. Since both leading powers are smaller than $(0,n)$, that is, $(0,n-1),(m-1,0)\prec (0,n)$, then the only possibility is that $X_{\omega_{-1}}(f)$ and $X_{\omega_{0}}(f)$ are not reducible modulo $\{f\}$. In particular, $lp(X_{\omega_{-1}}(f))=lp(h_{-1})=(0,n-1)$ and $lp(X_{\omega_{0}}(f))=lp(h_{0})=(m-1,0)$. These last equalities imply that $X_{\omega_{-1}}(f)=h_{-1}$ and $X_{\omega_{0}}(f)=h_0$. Thus, by Equation \eqref{eq:ap:dx}
$$
-f_y=X_{dx}(f)=Ah_{-1}+Bh_0.
$$
This shows that $f_y\in (B)$. In the same way, we find that $f_x\in (B)$. We are left to show that $f\in (B)$. Since, we already know that $f_x,f_y\in (B)$, it is equivalent to show that 
$$
f-\left(\frac{1}{m}xf_x+\frac{1}{n}yf_y\right)=\sum_{\substack{\alpha,\beta\geq 0\\ n\alpha+m\beta >nm}}\tfrac{nm-n\alpha-m\beta}{nm}z_{\alpha\beta}x^\alpha y^\beta\in (B).
$$
We are going to define a sequence of functions $f_\ell,g_\ell$ for $\ell\geq 1$ satisfying the following two conditions:
\begin{itemize}
\item $f_\ell=f-g_\ell$, with $g_\ell\in (B)$.
\item $g_\ell=X_{\eta_\ell}(f)-pf$, for some $\eta_\ell\in\Omega_{\mathbb C^2,\mathbf 0}^1$ and $p\in\mathbb C\{x,y\}$.
\end{itemize}

We start with $g_1=\left(\frac{1}{m}xf_x+\frac{1}{n}yf_y\right)$ and $\eta_1=\frac{1}{m}xdy-\frac{1}{n}ydx$, that is, $f_1=f-g_1$.  Now, we proceed in a inductive way. There are three cases:
\begin{enumerate}
\item[a)] $f_\ell=0$.
\item[b)] $lp(f_\ell)$ is not divisible by the leading power of any element of $B$.
\item[c)] $lp(f_\ell)$ is divisible by the leading power of some element of $B$.
\end{enumerate}

\emph{Case} a): we are done, since we can write $f=g_\ell\in (B)$.

\bigskip

\emph{Case} b): write $lp(f_\ell)=(\alpha,\beta)$, which by assumption is not divisible by any leading power of any element of $B$.  In particular, $(\alpha,\beta)$ is not divisible by $lp(h_{-1})=(0,n-1)$. Hence, $(\alpha,\beta)$ is not divisible by $(0,n)$. This means that $-f_\ell$ is a final reduction of $X_{\eta_\ell}(f)$ modulo $\{f\}$. Therefore, by  Proposition \ref{prop:impli:multiplicidad}, we obtain that $\nu_C(\eta_\ell)=\lambda=n(\alpha+1)+m(\beta+1)-nm$. Since $\lambda$ is a differential value, we have that $\lambda=\lambda_i+np+mq$ for some $-1\leq i\leq s$ and $p,q\geq 0$. 

Set $(c_i,d_i)=lp(h_i)$, this means that $\lambda_i=n(c_i+1)+m(d_i+1)-nm$. Thus, we find that $n\alpha+m\beta=n(c_i+p)+m(d_i+q)$. Therefore,
$$
\alpha=c_i+p-km; \quad \beta=d_i+q+kn; \quad k\in\mathbb Z.
$$
If $k>0$, then $\beta>n$ and $(\alpha,\beta)$ is divisible by $(0,n-1)=lp(h_{-1})$ which is a contradiction. If $k<0$, then $(\alpha,\beta)$ is divisible by $(m-1,0)= lp(h_0)$ again a contradiction. Finally, if $k=0$, then $(\alpha,\beta)=(c_i+p,d_i+q)$ and $(\alpha,\beta)$ is divisible by $(c_i,d_i)$ which is another contradiction. Thus, Case b) cannot happen.

\bigskip

\emph{Case} c): again, set $(\alpha,\beta)=lp(f_\ell)$ and consider $(c_i,d_i)=lp(h_i)$, such that $(\alpha,\beta)$ is divisible by $(a_i,b_i)$. Then, we write 
$$
\eta_{\ell+1}=\eta_\ell+\mu x^{\alpha-c_i}y^{\beta-d_i}\omega_i,
$$
where $\mu$ is the unique constant such that  the function 
$$
g_{\ell+1}=g_\ell+\mu x^{\alpha-c_i}y^{\beta-d_i}h_i,
$$ 
satisfies that  $lp(f-g_{\ell+1})>lp(f-g_\ell).$ 

By hypothesis, we can write $g_\ell=X_{\eta_\ell}(f)-p'f$ with $p'\in\mathbb C\{x,y\}$. Additionally, by definition of $h_i$, we have that 
$$
\mu x^{\alpha-c_i}y^{\beta-d_i}h_i=X_{\mu x^{\alpha-c_i}y^{\beta-d_i}\omega_i}(f)-p''f,\quad p''\in\mathbb C\{x,y\}. 
$$
Thus, if we write $p=p'+p''$, we get that $g_{\ell+1}=X_{\eta_{\ell+1}}(f)-pf$. Finally, we put $f_{\ell+1}=f-g_{\ell+1}$ and we restart the process.

This entire procedure shows that a final reduction of $f_1$ modulo $B$ is 0, implying that $f_1\in (B)$, in particular, we find that $f\in (B)$. This concludes the proof of Statement 1.

\bigskip

\emph{Statement} 2: Take $g\in \mathcal J(f)$ such that $lp(g)=(\alpha,\beta)$ is not divisible by any leading power of any element of $B$. By Statement 1, we can write 
$$
g=\sum_{i=-1}^s g_i h_i,
$$
for some functions $g_i$ with $i=-1,0,\ldots,s$. We consider the 1-form $\omega=\sum_{i=-1}^sg_i\omega_i$. As in Statement 1, the assumption of non divisibility implies, in particular, that $(\alpha,\beta)$ is not divisible by $(0,n-1)=lp(h_{-1})$. Hence, $lp(g)$ is not divisible by $(0,n)$ and $g$ is a final reduction of $X_\omega(f)$.

By Proposition \ref{prop:impli:multiplicidad}, $\nu_C(\omega)=\lambda=n(\alpha+1)+m(\beta+1)-nm$. Repeating exactly the same arguments that in Case b) from Statement 1, we find a contradiction with the assumption that $(\alpha,\beta)$ is not divisible by any leading power of the elements of $B$.

\bigskip

\emph{Statement} 3: Consider $h_i,h_j\in B$, with $i\neq j$. Assume that $(c_j,d_j)=lp(h_j)$ is divisible by $lp(h_i)=(c_i,d_i)$. Since $h_i,h_j$ are, respectively, final reductions of $X_{\omega_i}(f)$ and $X_{\omega_j}(f)$ modulo $\{f\}$, then by Proposition \ref{prop:impli:multiplicidad}, we have that 
\begin{eqnarray*}
\nu_C(\omega_i)&=&\lambda_i=n(c_i+1)+m(d_i+1)-nm\\
\nu_C(\omega_j)&=&\lambda_j=n(c_j+1)+m(d_j+1)-nm
\end{eqnarray*}
Since $(c_i,d_i)\mid (c_j,d_j)$, we have that $c_j-c_i\geq 0$ and $d_j-d_i\geq 0$. Therefore,   $\lambda_j-\lambda_i=n(c_j-c_i)+m(d_i-d_j)\in\Gamma_C$. This contradicts the fact that $\lambda_i$ and $\lambda_j$ are two different elements of the basis of $\Lambda_C$.
\end{proof}

\Addresses
\end{document}